\documentclass[12pt]{amsart}

\usepackage{amsmath,amssymb}

\newtheorem{thm}{Theorem}
\newtheorem{lem}[thm]{Lemma}
\newtheorem{prop}[thm]{Proposition}
\newtheorem{cor}[thm]{Corollary}
\theoremstyle{definition}
\newtheorem{rem}{Remark}

\newcommand{\Z}{{\mathbf Z}}
\newcommand{\Q}{{\mathbf Q}}

\newcommand{\R}{{\mathbf R}}
\newcommand{\C}{{\mathbf C}}
\newcommand{\E}{{\mathbb E}}

\newcommand{\vol}{\operatorname{vol}}

\newcommand{\SIGMA}{\sigma}
\newcommand{\Sp}{ \Omega_p}
\newcommand{\Sa}{ \Omega_{q_1}}
\newcommand{\Sb}{ \Omega_{q_2}}

\newcommand{\Sq}{ \Omega_q}

\newcommand{\Zq}{ \Z/q\Z}

\newcommand{\Za}{ \Z/q_1\Z}
\newcommand{\Zb}{ \Z/q_2\Z}

\newcommand{\vh}{{\bf h}}

\newcommand{\pr}{\text{Prob}}
\newcommand{\iif}{\quad\text{if}\ }
\newcommand{\om}{\Omega}
\newcommand{\vo}{\vol}

\newcommand{\eps}{\varepsilon}

\newcommand{\fp}{\mathbb{F}_p}
\newcommand{\kk}{\fp}
\newcommand{\kbar}{\overline{\mathbb{F}_p}}
\newcommand{\gal}{\operatorname{Gal}}
\newcommand{\fix}{\operatorname{Fix}}
\newcommand{\frob}{\operatorname{Frob}}
\newcommand{\sgn}{\operatorname{sgn}}
\newcommand{\m}{\mathfrak{m}}
\newcommand{\M}{\mathfrak{M}}

\newcommand{\jacobi}[2]{\genfrac{(}{)}{}{}{#1}{#2}}

\title{Poisson statistics via the Chinese remainder theorem}

\author{Andrew Granville}
\address{
D\'epartement de Math\'ematiques et statistique\\
                           Universit\'e de Montr\'eal\\
                           CP 6128 succ Centre-Ville\\
                           Montr\'eal QC  H3C 3J7\\
                           Canada}
\email{andrew@DMS.UMontreal.CA}

\author{P\"ar Kurlberg}
\address{Department of Mathematics\\
KTH\\
SE-100 44 Stockholm\\
Sweden}
\email{kurlberg@math.kth.se}

\thanks{A.G. has been supported in part by the National Science
  Foundation and by NSERC (Canada) during the preparation of this
  paper.  P.K. supported in part by 
the National Science Foundation,
the Royal Swedish Academy of Sciences, 
  and the Swedish Research Council.}

\subjclass{Primary 11N69, Secondary 11K36}

\begin{document}

\begin{abstract}
  We consider the distribution of spacings between consecutive
  elements in subsets of $\Z/q\Z$ where $q$ is highly composite and
  the subsets are defined via the Chinese remainder theorem.  We give
  a sufficient criterion for the spacing distribution to be Poissonian
  as the number of prime factors of $q$ tends to infinity, and as an
  application we show that the value set of a generic polynomial
  modulo $q$ have Poisson spacings.  We also study the spacings of
  subsets of $\Z/q_1q_2\Z$ that are created via the Chinese remainder
  theorem from subsets of $\Z/q_1\Z$ and $\Z/q_2\Z$ (for $q_1,q_2$
  coprime), and give criteria for when the spacings modulo $q_1q_2$
  are Poisson.  Moreover, we also give some examples when the spacings
  modulo $q_1q_2$ are not Poisson, even though the spacings modulo
  $q_1$ and modulo $q_2$ are both Poisson.

\end{abstract}

\date{February 10, 2005}
\maketitle

\section{Introduction }

Let $1=x_1<x_2<\dots < x_m<q$ be the set of
squares\footnote{An integer $x$ is a square mod $q$ if there exists
$y$ for which $y^2\equiv x \pmod
q$.}
modulo a 
large integer $q$. If $q=p$ is prime then $m=(p-1)/2$; that is,
roughly half of the integers mod $p$ are squares, so an integer
chosen at random is square with probability close to $1/2$.   So
do the squares appear as if they are ``randomly distributed'' (if
one can appropriately formulate this question)?  For instance, if
one chooses a random square $x_i$ mod $p$, what is the probability
that $x_{i+1}-x_i=1$, or $2$, or $3,\dots$? Is it the same as for
a random  subset of the integers?  In 1931 Davenport
\cite{Davenport-quadratic} showed that the answer is ``yes'' by
proving that the probability that $x_{i+1}-x_i=d$ is $
1/2^d+o_p(1)$. (Note that if one takes a random subset $S$ of
$[1,n]$ of size $n/2$ then the the proportion of $x\in S$ such
that the next smallest element of $S$ is $x+d$, is
 $\sim 1/2^d$ with probability 1.)

If $q$ is odd  with $k$ distinct prime factors, then
$m=\phi(q)/2^k$. The average gap, $s_q$, between these squares is
now a little larger than $2^k$, which is large if $k$ is large; so
we might expect that the probability that $x_{i+1}-x_i=1$ becomes
vanishingly small as $k$ gets larger. Hence, to test whether the
squares appear to be ``randomly distributed'', it is more
appropriate to consider $(x_{i+1}-x_i)/s_q$. If we have $m$
integers randomly chosen from $1,2,\dots ,q-1$ then we expect that
the probability that $(x_{i+1}-x_i)/s_q>t$  is $\sim e^{-t}$ as
$q, s_q \to \infty$. In 1999/2000 Kurlberg and Rudnick
\cite{squares1,squares2} proved that this is true for the squares
mod $q$.

To a number theorist this is reminiscent of Hooley's 1965 result
\cite{Hooley2,Hooley3} in which he proved that the set of integers
coprime to $q$ appear to be ``randomly distributed'' in the same
sense, as the average gap $s_q=q/\phi(q)$ gets large\footnote{
Under a similar assumption, namely that $s_p = (p-1)/\phi(p-1)$ tends to
infinity, Cobeli and Zaharescu has shown \cite{cobeli-zaharescu} that
the spacings 
between primitive roots modulo $p$ becomes Poissonian as $p$ tends to
infinity along primes.}.

In both of these examples the sets of integers $\Sq \subset \Zq$
are obtained from sets of integers $\Omega_{p^a} \subset \Z/p^e\Z$
(for each prime power $p^e\| q$) by the Chinese Remainder Theorem
(that is $a\in \Sq$ if and only if $a\in \Omega_{p^e}$ for all
$p^e\| q$). We thus ask whether, in general, sets $\Sq \subset
\Zq$ created from sets $\Omega_{p^a} \subset \Z/p^e\Z$ (for each
prime power $p^e\| q$) by the Chinese Remainder Theorem appear (in
the above sense) to be ``randomly distributed, at least under some
reasonable hypotheses? This question is inspired by the  Central
Limit Theorem which tells us that, incredibly, if we add enough
reasonable probability distributions together then we obtain a
generic ``random'' distribution, such as the Poisson or Normal
distribution.

Let us be more precise. For simplicity we restrict our attention to
squarefree $q$. Suppose that for each prime $p$ we are given a subset
$\Sp \subset \Z/p\Z$. For $q$ a squarefree integer, we define $\Sq
\subset \Z/q\Z$ using the Chinese remainder theorem; in other words,
$x \in \Sq$ if and only if $x \in \Sp$ for all primes $p$ dividing
$q$. Let $s_q = q/|\Sq|$ be the average spacing between elements of
$\Sq$, and $r_q=1/s_q=|\Sq|/q$ be the probability that a randomly
chosen integer belongs to $\Sq$. Let $1=x_1<x_2<\dots < x_m<q$ be the
elements of $\Sq$, and define $\Delta_j=(x_{j+1}-x_j)/s_q$ for all
$1\leq j\leq m-1$.  For any given real numbers $t_1,
t_2,\dots, t_k\geq 0$ define $\pr_q(t_1,\dots, t_k)$ to be the
proportion of these integers $j$ for which $\Delta_{j+i}>t_i$ for each
$i=1,2,\dots, k$.\footnote{By letting $x_j = x_{j \pmod m}$ and
  $\Delta_j = \Delta_{j \pmod m}$ for any $j \in \Z$ we obtain the
  distribution of spacings ``with wraparound'', but in the limit
  $|\Sq| \to \infty$, $\pr_q(t_1,\dots, t_k)$ is independent of
  whether spacings are considered with or without wraparound.}

Suppose that $Q$ is an infinite set of squarefree, positive
integers which can be ordered in such a way that $s_q\to \infty$.
We say that the spacings between elements in the sets $\Sq$ for
$q\in Q$ become {\sl Poisson distributed}  if for any $t_1,
t_2,\dots, t_m\geq 0$
$$
\pr_q(t_1,t_2, \ldots, t_m) \to  e^{-(t_1+t_2+ \ldots + t_m)} \
\text{\rm as} \ s_q\to \infty, \ q\in Q .
$$

For a given vector of integers $\vh=(h_1,h_2,\dots,h_{k-1})$, let
$h_0=0$ and define the counting function\footnote{The counting
  function is defined for $\vh$ modulo $q$, so implicitly we consider
  gaps with wraparound.} for $k$-tuples mod $q$ 
as
$$
N_k(\vh,\Sq) = \#\{ t\pmod q: \ t+h_i \in \Sq \ \text{\rm for} \
0\leq i\leq k-1 \} .
$$
Note that the average of  $N_k(\vh,\Sq)$ (over all possible $\vh$)
is $r_q^k\ q$.

Our main result shows that if for each fixed $k$, the $k$-tuples
of elements of $\Omega_p$ are well-distributed for all
sufficiently large primes $p$, then indeed the sets $\Sq$ become
Poisson distributed.

\begin{thm}
\label{thm:main} Suppose that we are given subsets $\Sp \subset
\Z/p\Z$ for each prime $p$. For each  integer $k$, assume that
\begin{equation}
  \label{eq:0}
N_k(\vh,\Sp) = r_p^k \cdot  p\ (1+ O_k( (1-r_p)p^{-\epsilon}))
\end{equation}
provided that  $0, h_1, h_2, \ldots, h_{k-1}$ are distinct mod $p$.
If $s_p = p^{o(1)}$ for all primes $p$, then the spacings
between elements in the sets $\Sq$  become Poisson distributed  as
$s_q\to \infty$.
\end{thm}
\begin{rem}
Theorem~\ref{thm:main-bis} in section~4 actually gives something a
little more explicit and stronger.
\end{rem}
 From the theorem, we easily recover the result of Hooley, since
for $\Sp=\{ 1,2,\dots ,p-1\}$  we have $r_p=1-1/p$ and thus
$$
N_k(\vh,\Sp) = p-k =  r_p^k \cdot p \left(1+ O_k\left( \frac
{1-r_p}p\right)\right);
$$
and a generalization of the result Kurlberg-Rudnick using Weil's
bounds for the number of points on curves:

\begin{cor}
\label{cor:powers}
Fix an integer $d$ and let $\Sq$ be the set of $d$-th powers
modulo $q$. Then the spacings between elements in the sets $\Sq$
become Poisson distributed  as $s_q\to \infty$.
\end{cor}

Another situation where we may apply Weil's bounds is to the sets

\noindent $\{ x \mod q:\ \text{There exists } y \mod q \
\text{such that } y^2\equiv x^3+ax+b \pmod q\}$ for any given
integers $a,b$; and indeed to coordinates of any given
non-singular  hyperelliptic curve.  Thus we may deduce the analogy to
Corollary~\ref{cor:powers} in these cases.

 In section~\ref{sec:general-polynomials}
we also show that the spacings between residues mod $q$ in the
image of a polynomial having $n-1$ distinct critical values\footnote{The
  critical values 
 of $f$ is the set $\{f(\xi) : \xi \in \C, f'(\xi) = 0 \}$.}
(a generic condition) become Poisson distributed  as $s_q\to
\infty$:
\begin{thm}
\label{thm:gaps-of-polynomials}
  Let $f$ be a polynomial of degree $n$ with integer coefficients.
 Regarding $f$ as
  a map from $\Zq$ into itself, define $\Sq$ to be the image of $f$
  modulo $q$, i.e., $\Sq := \{ x
 \pmod q : \text{ there exists } y \pmod q \text{ such that } f(y)
 \equiv x \pmod q\}$.
If $f$ has $n-1$ distinct critical values, then the
 spacings between elements in the sets $\Sq$  become Poisson
 distributed  as $s_q\to \infty$.
\end{thm}

\begin{rem}
Theorem~\ref{thm:gaps-of-polynomials} is
  true for all polynomials, 
  but the proof of this is considerably more complicated  and
will appear in a separate paper. 
In fact, there are polynomials for which (\ref{eq:0})
does not hold - see section 4.2 for more details.
We also note that if $f$ has $n-1$ distinct critical values, 
Birch and Swinnerton-Dyer have proved \cite{bsd-valueset}
that 
$$|\Sp| = |\{x \in \fp : x = f(y) \text{ for some $y \in \fp$ }\}|
= c_n p + O_n(p^{1/2})
$$
where 
$$
c_n = 1 - \frac{1}{2} + \frac{1}{3!} - \ldots - (-1)^n \frac{1}{n!}
$$
is the truncated Taylor series for $1-e^{-1}$.  
(Note that $n! \cdot (1-c_n)$ is the "nth derangement number" from
combinatorics, so $c_n$ can be interpreted as the probability that a
random permutation $\sigma \in S_n$ has at least one fixed point.  In fact,
this is no coincidence - for these polynomials the Galois group of
$f(x)-t$, over $\fp(t)$, equals $S_n$, and the proportion of elements
in the image of $f$, up to an error $O(p^{-1/2})$, equals the
proportion of elements in the Galois group fixing  at least one root.)
Since the expected
cardinality of the image of a random map from $\fp$ to $\fp$ is
$p \cdot (1-e^{-1})$, the above result can be interpreted as saying
that the cardinality of the image of a generic polynomial (of large
degree) behaves as that of a random map.  Their result also implies that
$s_q \to \infty$ as the number of prime 
factors of $q$ tends to infinity.
\end{rem}

In Theorem~\ref{thm:main} we  proved that if {\sl all} $k$-tuples in
$\Sp$ are 
``well-distributed'' (in the sense of (1)) for all primes $p$ then
the $\Sq$ become Poisson distributed as $s_q\to \infty$. Perhaps
though one needs to make less assumption on the sets $\Sp$? For
example, perhaps it suffices to simply assume an averaged form of
(1), like
$$
\frac 1{p^{k-1}} \sum_{\vh} \left| \frac{N_k(\vh,\Sp)}{r_p^k\ p} -
1 \right| \ll_k  (1-r_p)p^{-\epsilon}
$$
where the sum is over all $\vh$ for which $0, h_1, h_2, \ldots
h_{k-1}$ are distinct mod $p$. We have been unable to prove this
as yet.

In the central limit theorem, where one adds together lots of
distributions to obtain a normal distribution, the hypotheses for the
distributions which are summed is very weak. So perhaps in our problem
we do not need to make an assumption which is as strong 
as (1)?  In section~\ref{crt-poisson-general} we suppose that we are
given sets $\Sa$ and $\Sb$ of residues modulo $q_1$ and $q_2$ (with
$(q_1,q_2)=1$), and try to determine whether the spacings in $\Sq$
(where $q=q_1q_2$) is close to a Poisson distribution. We show that
under certain natural hypotheses the answer is ``yes''.
These take the form: If $\Omega_{q_1}$ is suitably "strongly Poisson"
then $\Omega_q$ is Poisson if and only if $\Omega_{q_2}$ is Poisson
with an appropriate parameter.

On the other hand,  if we allow the sets to be correlated, then
the answer can be ``no''. In section~\ref{sec:counterexamples} we
give three examples in which the distribution of points in $\Sq$
is not consistent with that of a Poisson distribution.  The
constructions can be roughly described as follows:
\begin{itemize}
\item  $\Sa$ is random and small, and   $\Sb=\{ a:\ 1\leq a\leq q_2/2\}$ .

\item $\Sb=\Sa$ is a random subset of $\{ 1,2,\dots ,q_1\}$ where $q_2=q_1+1$.

\item Each $\Omega_{q_i}$ is a random subset of $\{ a: 1\leq a\leq q_i,\
m|a\}$ for $i=1,2$, with integer $m\geq 2$.
\end{itemize}

\subsection{Acknowledgments}
P.K. would like to thank J.~Brzezi\'nski, T.~Ekedahl, M.~Jarden, and
Z.~Rudnick for  helpful discussions.

\section{Poisson statistics primer}
\label{sec:poisson-statistics-primer}

Given a positive integer $q$ and a subset  $\Sq \subset \Zq$, let
$s_q  = q/|\Sq|$ be the average gap between consecutive elements
in $\Sq$. One can view $r_q=1/s_q$ as  the probability that  a
randomly selected element in $\Z/q\Z$ belongs to $\Sq$.

If $0< x_1<x_2<\dots $ are the positive integers belonging to
$\Sq$ then define $\Delta_j=(x_{j+1}-x_j)/s_q$ for all $j\geq 1$;
we are interested in the statistical behavior of these gaps as $q
\to \infty$, along some subsequence of square free integers. We
define the (normalized) {\em
  limiting spacing distribution}, if it exists, as a probability
measure $\mu$ such that
$$
\lim_{q \to\infty} \frac{\#\{j :\ 1\leq j\leq |\Sq|,\ \Delta_j \in
I \}}{|\Sq|} = \int_{I} d \mu(x)
$$
for all compact intervals $I \subset \R^+$. If $d \mu(x) = e^{-x}
\, dx$ and the gaps are independent
(i.e.,  that $k$ consecutive gaps are independent for any $k$),
the limiting spacing 
distribution is said to be {\sl Poissonian}.
This can be characterized (under fairly general conditions) as
follows: For any fixed $\lambda>0$ and integer $k\geq 0$, the
probability that there are exactly $k$ (renormalized) points in a
randomly chosen interval of length $\lambda$, is given by
$\frac{\lambda^k e^{-\lambda}}{k!}$ (see
\cite{billingsley-probability-and-measure}, section 23.)

We shall use a characterization of the Poisson distribution  which
is relatively easy to work with:\ The {\sl $k$-level correlation}
for a compact set $X \subset \{ x \in \R^{k-1} : 0<x_1 < x_2 <
\ldots < x_{k-1}\}$ is defined as
\begin{equation}
  \label{eq:6}
R_k(X, \Sq) = \frac{1}{|\Sq|} \sum_{\vh \in  s_qX \cap \Z^{k-1} }
N_k(\vh, \Sq) .
\end{equation}
Note that we ensure that $0<h_1<\ldots <h_{k-1}$ else $N_k(\vh,
\Sq)=N_\ell(\vh', \Sq)$ where $0<h'_1<\ldots <h'_{\ell-1}$ are the
distinct integers amongst $0,h_1,\dots, h_{k-1}$.

Now for any positive real numbers $b_1, b_2,\dots , b_{k-1}$
define
$$
B(b_1, b_2,\dots , b_{k-1}):= \{ x \in \R^{k-1} :
0<x_i-x_{i-1}\leq b_i \ \text{for} \ i=1,2,\dots k-1\}
$$
where we let $x_0 = 0$.
Let $\mathbb B_k$ be the set of such (not necessarily rectangular)
boxes.

Suppose we are given a sequence of integers $Q=\{
q_1,q_2,\dots \}$ with $s_{q_i}\to \infty$ as $i\to \infty$. Then
(e.g., see Appendix~A of \cite{squares1}) 
the spacings of the elements in $\Omega_{q_n}$ become
Poisson  as
$n\to \infty$ if and only if  for each integer $k\geq 2$ and box
$X\in \mathbb B_k$,
$$
R_k(X, \Omega_{q_n})\to \vol(X) \  \ \text{\rm  as} \ n\to \infty
.
$$

It will be useful to include a further definition along similar
lines. Suppose $\theta_n$ is a positive real number for each $n$.
We say that the spacings of the elements in $\Omega_{q_n}$ become
Poisson with
parameter $\theta_n$ as $n\to \infty$ if and only if 
for each integer $k\geq 2$ and box $X\in \mathbb B_k$,
$$
R_k(\theta_nX, \Omega_{q_n})\to \vol(\theta_nX) \  \ \text{\rm
as} \ n\to \infty .
$$
Notice that ``Poisson with parameter $1$''  is the same thing as
``Poisson''.  (In fact, Poisson with any bounded parameter is the same
as Poisson.)

\subsection{Correlations for randomly selected sets}
\label{sec:corr-rand-select} Let $X_1,\ X_2,\dots , X_q$ be
independent Bernoulli random variables with parameter $1/\sigma
\in (0,1)$. In other words, $X_i = 1$ with probability $1/\sigma$,
and $X_i = 0$ with probability $1-1/\sigma$.  Given an outcome of
$X_1, X_2, \ldots, X_q$, we define $\Sq \subset \Zq$ by letting $i
\in \Sq$ if and only if $X_i=1$.  Note that the expected
average gap
is then given by $\sigma$.
Below we write $R_k(x,q)$ for $R_k(x,\Sq)$.
\begin{lem}
\label{lem:random-correlations} As we vary over all subsets of
$\Zq$ with the probability space as above, we have
$$
\E(R_k(X,q)) =\vo(X)+O_k\left(1/\sigma + \sigma/q \right)
$$
and
$$
\E\left(\bigl(R_k(X,q)-\vo(X)\bigr)^2\right)\ll_k 1/\sigma +
\sigma/q
$$
\end{lem}

\begin{proof}
Using conditional expectations we write
\begin{multline*}
\E(R_k(X,q))
=\sum_{r = k}^q\pr(|\Sq|=r)\ \E(R_k(X,q):|\Sq|=r) \\
=\sum_{h\in \sigma X \cap \Z^{k-1}}\ \sum_{r = k}^q\ \frac{\pr(|\Sq|=r)}{r}\
\sum_{i=1}^q\ \E\left(x_i x_{i+h_1}\dots x_{i+h_{k-1}}:|\Sq|=r\right)
\end{multline*}
Now, the number of ways to have
$|\Sq|=r$ is $\binom qr$, and the number of ways to
have $|\Sq|=r$
with $i, i+h_1,\dots, i+h_{k-1}\in\Sq$ is $\binom{q-k}{r-k}$. Therefore,
$$
\E\left(x_i x_{i+h_1}\dots x_{i+h_k}:|\Sq|=r\right)
=\dbinom{q-k}{r-k}\Big/ \dbinom qr
$$
Note that $R_k(X,q)=0$ if $|\Sq|\le k-1$, and
$$
\pr (|\Sq|=r)=\binom qr  \sigma^{-r}(1-1/\sigma)^{q-r}.
$$
Taking $q\geq 4k$ with $q/\sigma$ large, we obtain
\begin{multline*}
\E(R_k(X,q))
=\left| \sigma X \cap \Z^{k-1}\right|
\sum_{r =  k}^q\ \frac 1r\  \sigma^{-r}(1-1/ \sigma)^{q-r} 
q\cdot\binom{q-k}{r-k}\\
=q \sigma^{-k}\left( \sigma^{k-1}\ \vo(X)+O( \sigma^{k-2})\right)
\cdot\sum_{r = k}^q\ \frac {\sigma^{k-r}}{r}\
 (1-1/ \sigma)^{(q-k)-(r-k)}\binom{q-k}{r-k}
\\
=(q/ \sigma) \left( \vol(X)+O( 1/\sigma) \right) \cdot\sum_{R =
0}^Q\ \frac 1{R+k}\   (1/\sigma)^R (1-1/ \sigma)^{Q-R}\binom{Q}{R}
\end{multline*}
where $Q=q-k$ and $R=r-k$. Now
$$
\frac 1{R+k}=\frac 1{R+1}+O\left(\frac k{(R+1)(R+2)}\right) ,
$$
so the last sum is
$$
  \frac  \sigma{(Q+1)}\left(1-(1-1/ \sigma)^{Q+1}\right)
+O\left(\frac {k  \sigma^2}{Q^2}\right) =\frac
\sigma{q}\left(1+O_k\left(\frac  \sigma{q}\right)\right),
$$
since $(Q/\sigma)^A (1-1/\sigma)^Q   \ll_A 1$, and thus
$$
\E(R_k(X,q))
=\vo(X)+O\left(1/ \sigma+\sigma/q \right).
$$
For the variance,  note that 
\begin{multline*}
\E\left(R_k(X,q)^2\right)
=
\sum_{r= k}^q \pr \left(|\Sq| = r\right)\E\left(R_k(X,q)^2:|\Sq|=r\right)
\end{multline*}
\begin{multline*}
=\sum_{r= k}^q \binom qr  \sigma^{-r}(1-1/ \sigma)^{q-r}
\frac1{r^2}
\cdot\\
\cdot
\sum_{h, H\in  \sigma X \cap \Z^{k-1}} \
\sum_{i,j=1}^q
\E\left(x_i x_{i+h_1} x_{i+h_2}\dots x_{i+h_{k-1}}
x_j x_{j+H_1}\dots x_{j+H_{k-1}}:|\Sq|=r\right)
\end{multline*}
If there are $l$ distinct elements in $\{i, i+h_1,\dots, h_{k-1}, j,
j+H_1, \ldots,  j+H_{k-1}\}$
then the expectation is
$$
\binom{q-l}{r-l}\biggl/\binom{q}{r}.
$$

Given $\alpha,\beta,h$ and $H$ there is a solution to $i+h_\alpha=j+H_\beta$ for
$O(k^2q)$ values of
$i$ and $j$.
Thus our main term is
$$
\left(q^2+O_k(q)\right)\binom{q-2k}{r-2k}\biggl/\binom{q}{r}.
$$
We treat the other terms as follows:  Fix $d$ and consider $i$ and $j$
with $j\equiv i+d\pmod q$.
Select $u_1,\dots, u_m, v_1,\dots,v_m$ with $h_{u_t}\equiv H_{v_t}+d\pmod q$.
The number of choices for $i$ and $j$ is $q$.  $H$ can be chosen
freely and so can $k-m-1$ of the coordinates of $h$. The total number of
choices is thus
$$
\asymp_{X,k}q  \sigma^{k-1} \sigma^{k-m-1}
$$
Moreover the number of choices for $d$ is $\asymp_X  \sigma$.
Therefore, since $l=2k-m$, we have\footnote{We use the convention that
$\binom{n}{k}=0$ if $k<0$.}
\begin{multline*}
\E\left(R_k(X,q)^2\right)
=
\sum_{r =  k}^q
\frac{ \sigma^{-r}(1-1/ \sigma)^{q-r}}{r^2} 
\times \\ \times 
\Bigg(
\left| \sigma X\cap\Z^{k-1}\right|^2 \binom{q-2k}{r-2k} \left(q^2+O(q)\right)
+ 
O\left( \sum_{m=1}^k\ \binom{q-2k+m}{r-2k+m}q  \sigma^{2k-1-m}\right) \Bigg)\\
= \left(q^2+O(q)\right)\left( \sigma^{k-1}\vo(X)+O(
\sigma^{k-2})\right)^2\
\sum_{r=2k}^q\binom{q-2k}{r-2k}\frac1{r^2}\ \sigma^{-r}(1-1/
\sigma)^{q-r}
\\
+ O\left( \sum_{m=1}^k q \sigma^{2k-1-m} \sum_{r =  2k-m}^q
\binom{q-2k+m}{r-2k+m} \frac{ \sigma^{-r}(1-1/ \sigma)^{q-r}}{r^2}
\right) .
\end{multline*}
Now, for $k\leq \ell \leq 2k$ take $Q=q-\ell$ and $R=r-\ell$, and
note that
$$
\frac 1{(R+\ell)^2}=\frac 1{(R+1)(R+2)}+O_k\left(\frac
1{(R+1)(R+2)(R+3)}\right) ,
$$
to obtain
\begin{multline*}
\sum_{r=\ell}^q\binom{q-\ell}{r-\ell}\frac1{r^2}\ \sigma^{-r}(1-1/
\sigma)^{q-r} 
\\
= \sigma^{-\ell}\sum_{R=0}^Q\binom{Q}{R}\frac
1{(R+\ell)^2}\
 (1/\sigma)^R(1-1/ \sigma)^{Q-R} 
 \sigma^{-\ell}
\cdot \\ \cdot
\left( \frac{ \sigma^2}{(Q+1)(Q+2)} + O_k \left( \frac{
\sigma^3}{q^3} \right) \right) 
\\= \frac{ \sigma^{2+2k-\ell}}{
\sigma^{2k}q^2} \left( 1+ O_k\left( \frac{ \sigma}{q} \right)
\right)
\end{multline*}
Substituting this in above gives
$$
\E\left(R_k(X,q)^2\right)
=\vo(X)^2+O\left(1/ \sigma+\sigma/q\right),
$$
and hence
\begin{multline*}
\E\left(\bigl(R_k(X,q)-\vo(X)\bigr)^2\right)
=
\E\left(\bigl(R_k(X,q)\bigr)^2\right)-\vo(X)^2
\\ = O\left( 1/ \sigma+\sigma/q \right).
\end{multline*}

\end{proof}

One can interpret this result as saying that almost all sets have
Poisson spacings.

\section{Correlations  via the Chinese remainder theorem}

\subsection{Counting solutions to congruences}
\label{sec:higher-correlations-preliminaries}

Suppose that $\Gamma=\{\gamma_{i,j}: 0\leq i\neq j \leq k-1
\mbox{  with } \gamma_{i,j}=\gamma_{j,i}\}$ is a given set of positive
  squarefree integers for which
\begin{equation}
\label{eq:y-star}
\gcd (\gamma_{i,j}, \gamma_{j,l}) \mbox{ divides } \gamma_{i,l}
\mbox{ for any distinct } i,j, l 
\end{equation}
Define
$$
\gamma_j :=
\underset{{0\leq i\leq j-1}}{\operatorname{LCM}}
\gamma_{i,j}
$$
and let
$$
\gamma(\Gamma) :=\gamma_1 \ldots \gamma_{k-1}
$$

Once one understands all this terminology one easily sees that
\begin{lem}
If $\sigma$ is a permutation of $\{1,\ldots, k-1\}$ and $\sigma(0)=0$
define $\gamma^{(\sigma)}_{i,j}=\gamma_{\sigma(i), \sigma(j)}$. Then
$\gamma^{(\sigma)} (\Gamma)=\gamma (\Gamma)$.
\end{lem}

Define $c(\Gamma)$ to be the squarefree product of the primes
dividing $\gamma(\Gamma)$, so that $c(\Gamma)$ divides $\gamma
(\Gamma)$, which divides $c(\Gamma)^{k-1}$.

Given a squarefree positive integer $c$, and a set of distinct
non-negative integers $h_0=0, h_1, h_2, \ldots, h_{k-1}$ 
let $\vh = (h_1, \ldots, h_{k-1})$ and define
\[
\gamma_{i,j}(\vh): =
\gcd (c, h_j-h_i) \mbox{ for } 0\leq i\neq j \leq k-1,
\]
and then  $\Gamma(\vh)$ accordingly.

For a given set $\Gamma$ and integer $c=c(\Gamma)$ define
\begin{multline}
M_\Gamma (H) := \#  \{ (h_0=0, h_1,\ldots, h_{k-1}) \in {\mathbb Z^k}:
\\
h_i \neq h_j \text{ for $i \neq j$, } 0 \leq h_i
\leq H \mbox{ for all } 0\leq i\leq k-1 \mbox{ and }
\Gamma(\vh)=\Gamma \}
\end{multline}

Finally for given integers $\gamma$ and $c$, with $c|\gamma|c^{k-1}$,
define
\begin{equation}
\label{eq:y-star-zero}
M_\gamma (H) :=
\sum_{\Gamma: \gamma(\Gamma)=\gamma}M_\Gamma (H).
\end{equation}

We wish to give good upper bounds of $M_\gamma (H)$. First note
that if $\gamma_{i,j}>H$, then $M_\Gamma (H)=0$ else
$\gamma_{i,j}| h_i-h_j$ and so $H<\gamma_{i,j}\leq |h_i-h_j|\leq
H$. Thus if $\gamma>H^{\binom{k}{2}}$ then $M_\gamma(H)=0$ else
$\max \gamma_{i,j}\geq \gamma^{1/\binom{k}{2}}>H$.

The Stirling number of the second kind, $S(k,\ell)$, is defined to
be the number of ways  of partitioning a $k$ element set into
$\ell$ non-empty subsets, and may be evaluated as
$$
S(k,\ell) = \frac{1}{(\ell-1)!}  \sum_{j=1}^\ell (-1)^{\ell-j}
\binom{\ell-1}{j-1} j^{k-1} .
$$
One can show  that $S(k,k-e)\leq \binom k2^e$.

\begin{lem}
\label{lem:gamma-bound}
$\# \{\Gamma: \gamma(\Gamma)=\gamma\} \leq \prod_{p^e\| \gamma}
S(k,k-e) \leq \binom k2 ^{\# \{p^e:p^e|\gamma\}}$.
\end{lem}

\begin{proof}
 For each prime $p$ dividing $\gamma$, we partition $\{0,
\ldots, k-1\}$ into subsets where $i$ and $j$ are in the same subset
if $p|\gamma_{i,j}$ (by (\ref{eq:y-star}) this is consistent). The bound
follows.
\end{proof}

Now we wish to bound $M_\Gamma (H)$.

\begin{prop}
\label{prop:m-gamma-bound}
We have
\[
M_\Gamma (H) \leq \prod^{k-1}_{i=1} \Big
  (\frac{H}{\gamma^{(\sigma)}_i} +1\Big ) \mbox{ for any } \sigma \in
  S_{k-1}.
\]
\end{prop}

\begin{proof}
Certainly we may rearrange the order, using $\sigma$, without changing
the question; so relabel $\sigma(i)$ as $i$. Now by induction on $k
\geq 1$, we have, for each given $(h_1,\ldots, h_{k-2}) \in
M_{\Gamma'}(H)$ where $\Gamma'$ is $\Gamma$ less all elements of the
form  $\gamma_{i,k-1}$ or $\gamma_{k-1,i}$ for $0 \leq i \leq k-1$,
that if $(h_1,\ldots, h_{k-1}) \in 
M_\Gamma (H)$, then 
$h_{k-1} \equiv h_i \mbox{ mod } \gamma_{i,k-1}$ for each $i$, $0\leq
i\leq k-2$ and so $h_{k-1}$ is determined modulo $\gamma_{k-1}$. Thus the
number of possibilities for $h_{k-1}$ is $\leq H/\gamma_{k-1}+1$, and
the result follows.
\end{proof}

\begin{cor}
\label{cor:m-gamma-bound-one}
We have
$$
M_\Gamma(H)
\leq
2^{k-1}H^{k-1}/\prod^k_{i=1} \min(\gamma_i,H)
$$
In particular,
\begin{equation}
  \label{eq:y-m-gamma}
M_\Gamma (H)\leq
\begin{cases}
  2^{k-1}H^{k-1}/\gamma & \text{if each $\gamma_i \leq H$} \\
  2^{k-1}H^{k-2} & \text{if any $\gamma_j \geq H$}
\end{cases}
\end{equation}
\end{cor}
{\em Remark:} 
When $k=2$  the first bound in  (\ref{eq:y-m-gamma}) is up to the
constant best possible. For $k=3$ things are immediately more
complicated. For suppose $\gamma_{0,1}, \gamma_{0,2},
\gamma_{1,2}$ are all coprime and each lies in the interval
$(T,2T)$ with $T>\sqrt{H}$. Then $\gamma_1 \approx T, \gamma_2>H$
and so $M_\Gamma(H) \leq 4H/T$ is what the corollary yields,
rather than what we might predict, $\approx H^2/T^3$. Thus this
``prediction'' cannot be true if $T>H^{2/3+\epsilon}$.

Next we look for a ``good'' re-ordering $\sigma$; select $\sigma(1)$
so as to maximize $\gamma_{\sigma(1),0}$. Now swap $\sigma(1)$
and 1 and then swap $\sigma(2)$ and 2 so as to maximize
$\operatorname{LCM}(\gamma_{\sigma(2),1}, \gamma_{\sigma(2),
  0})$. Proceeding  like this we obtain
\[
\gamma_r=LCM[\gamma_{r,0}, \gamma_{r,1}, \ldots , \gamma_{r,r-1}]\geq
LCM[\gamma_{j,0}, \gamma_{j,1}, \ldots \gamma_{j,r-1}] \mbox{ for all }
 j \geq r.
\]
Note that
\begin{equation}
  \label{eq:y-3}
\gamma_{r+1}
\leq
LCM[\gamma_{r,0}, \ldots, \gamma_{r,r-1}] \gamma_{r+1,r}
=
\gamma_r \gamma_{r+1,r} \leq H \gamma_r.
\end{equation}

Now in our general construction let   $I=\{i \in [1,\ldots, k-1]:
\gamma_i\leq H\}$ and write $D(\Gamma)=\prod^{k-1}_{i=1}
\min(\gamma_i,H)$ so that $M_\Gamma(H)\leq (2H)^{k-1}/D(\Gamma)$, and 
$D(\Gamma)=H^{k-|I|-1}D_I(\Gamma)$ where
$D_I(\Gamma)=\prod_{i\in I} \gamma_i$. Also, by (\ref{eq:y-3}) we have 
$\gamma_{r+1}\leq
H\gamma_r$, and thus
\[
\gamma=\gamma_1 \ldots \gamma_{k-1} \leq \prod_{i\in I} \gamma_i \cdot
\prod^{k-|I|-1}_{j=1} H^{1+j}=D_I(\Gamma) H^{\frac{1}{2} (k -
  |I|-1)(k-|I|+2)}.
\]
Let us suppose $|I|=\rho$ where $1\leq \rho \leq k-1$ (note that we
always have $\gamma_1 \leq H$). Then $1\leq D_I(\Gamma)\leq
H^\rho$. Write
$D_I(\Gamma) =H^{\rho \theta}$ for some $0\leq \theta \leq 1$. Thus
\begin{equation}
  \label{eq:y-1}
D(\Gamma)=H^{k-1-\rho+\rho \theta}
\end{equation}
and
\begin{equation}
  \label{eq:y-2}
\gamma
\leq H^{\rho \theta+\frac{1}{2} (k-\rho-1)(k-\rho+2)}
\leq H^{\frac{1}{2}(k-\rho-1)(k-\rho+2)+\rho}
\end{equation}

We note that $\frac{1}{2} (k-\rho-1)(k-\rho+2)+\rho$ is decreasing
in the range $1\leq \rho \leq k-1$.  Therefore if we choose $\tau$
in the range $1\leq \tau \leq k-1$ so that
\begin{equation}
  \label{eq:y-4}
H^{\frac{1}{2}(\tau-2)(\tau+1)+k+1-\tau}<\gamma \leq
H^{\frac{1}{2}(\tau-1)(\tau+2)+k-\tau}
\end{equation}
then $\rho\leq k-\tau$.

We wish to bound $D(\Gamma)$ from below.  By (\ref{eq:y-1}), we
immediately get
$$
D(\Gamma) \geq H^{k-1-\rho}
$$
Moreover, if for a given $\rho \leq k-\tau$, we have $\gamma \leq
H^{\frac{1}{2}(k-\rho-1)(k-\rho+2)+\rho \theta}$ then
$$
H^{\rho \theta} \geq \frac{\gamma}{H^{\frac{1}{2}(k-\rho-1)(k-\rho+2)}}
$$
and thus
$$
D(\Gamma) = H^{k-1-\rho} \cdot H^{\rho \theta} \geq \frac{\gamma
H^{k-1-\rho}}{H^{\frac{1}{2}(k-\rho-1)(k-\rho+2)}} =
\frac{\gamma}{H^{\frac{1}{2}(k-\rho-1)(k-\rho)}} .
$$

Since we are going to relinquish control of $\gamma$, other than
the size, we obtain the bound from the worst case. To facilitate
the calculation, we write $\gamma=H^\lambda, D(\Gamma)=H^\Delta$
and $\mu=k-1-\rho$ so that $k-2\geq \mu \geq \tau-1$.  With this
notation, (\ref{eq:y-4}) is equivalent to
\[
\frac{\tau^2}{2} -\frac{3\tau}{2} +k < \lambda \leq
\frac{\tau^2}{2} - \frac{\tau}{2} + k-1 .
\]
For a given $\lambda$ in our range we thus have, from the bounds
above,
\[
\Delta \geq \min_{\mu \geq \tau} \left( \max \left\{
\min_{\substack{\mu:\\
    \frac{1}{2} \mu (\mu+3)\geq \lambda}} \mu,
\min_{\substack{\mu:\\
    \frac{1}{2} \mu(\mu+3)\leq \lambda}} \lambda-\frac{1}{2} \mu
    (\mu+1)
\right\} \right) \geq u
\]
where we define $u$ to be the positive real number for which
\[
\frac{1}{2} u(u+3)=\lambda
\]
so that
$$
\left(u+\frac{3}{2}\right)^2=u(u+3)+\frac{9}{4} = 2\lambda +
\frac{9}{4} > \left(\tau-\frac{3}{2}\right)^2+2k \geq 2k+
\frac{1}{4},
$$
if $\tau$ is an integer. Note also that  $H^\Delta=D(\Gamma) \geq
H^{k-1-\rho} \geq H^{k-1-(k-\tau)}$ so that $\Delta \geq \tau-1$.
Therefore  $\Delta \geq \max (\tau-1,\sqrt{2k+1/4}-3/2)$. Thus we
have proved the following:

\begin{cor}
\label{cor:m-gamma-bound} Let $\tau$ be an integer $1 \leq \tau
\leq k$, and define $w(\tau)=\frac{1}{2} (\tau
-\frac{1}{2})^2+k-\frac{9}{8}$. If $H^{w(\tau-1)}<\gamma \leq
H^{w(\tau)}$ then
\[
M_\Gamma (H) \ll_k H^{k-\max\{\tau,\sqrt{2k+1/4}-1/2\}}.
\]
Note that $w(k-1)=k(k-1)/2$, and let
$\tau_1=[\sqrt{2k+1/4}-\frac{1}{2}]$. Combining this with
Lemma~\ref{lem:gamma-bound} and 
Corollary~\ref{cor:m-gamma-bound-one} gives that
\begin{multline*}
M_\gamma (H) \ll_k \prod_{p^e || \gamma} S(k,k-e) \cdot \\ \cdot 
\begin{cases}
H^{k-1}/\gamma &
\text{for $\hphantom{H^{w(0)}} \gamma \leq H,$} \\
 H^{k-2} &
\text{for $\hphantom{w(0)} H < \gamma \leq H^{\omega(0)}$} \\
H^{k+1/2-\sqrt{2k+1/4}}  &
\text{for $H^{w(0)}\,<\gamma \leq H^{w(\tau_1)}$} \\
H^{k-\tau} &
\text{for $H^{w(\tau-1)} < \gamma \leq H^{w(\tau)}$
$\tau_1+1 \leq \tau \leq k-1$}
\end{cases}
\end{multline*}
\end{cor}

\subsection{Proof of Theorem~\ref{thm:main}}

For $\vh \in \Z^{k-1}$, define the ``error term''
  $\eps_k(\vh,q)$ 
by
$$
N_k(\vh,q) = r_q^{k-1} |\Sq|(1 + \eps_k(\vh,q)) .
$$
We will need to use bounds on the size of $|\eps_k(\vh,p)|$, so
select $A_{p,k}$ so that
\begin{equation*}
|\eps_k(\vh,p)| \leq A_{p,k}
\end{equation*}
for all $\vh$ for which $0,h_1,\dots h_{k-1}$ are distinct mod
$p$. If $0,h_1,\dots h_{k-1}$ are not all distinct mod $p$ then
let  $\vh'$ be the set of distinct residues amongst $0, h_1,
\ldots, h_{k-1}$ mod $p$; if $\vh'$ contains $\ell \geq 1$
elements, then $N_k(\vh,p) = N_\ell (\vh',p)$ so that
\begin{equation}
  \label{eq:106}
\eps_k (\vh,p)=s^{k-\ell}_p-1 + s^{k-\ell}_p \eps_\ell(\vh',p).
\end{equation}
We will assume that $A_{p,k}$ is non-decreasing as $k$
increases\footnote{This is a benign assumption since we may
replace each $A_{p,k}$ by $\max_{\ell\leq k} A_{p,\ell}$.}.

For  $d>1$ a square free integer, put $e_k(\vh,1) = 1$ and
$$
e_k(\vh,d) = \prod_{p|d} \eps_k(\vh,p) ,
$$
so that
$$
N_k(\vh,q) = \prod_{p|q} r_p^{k-1} |\Sp| \left(1 + e_k(\vh,p)
\right) = r_q^{k-1} |\Sq| \sum_{d|q} e_k(\vh,d).
$$
With this notation
$$
R_k(X,\Sq) = \frac{1}{|\Sq|} \sum_{\vh \in s_q X \cap \Z^{k-1}}
N_k(\vh,q) = r_q^{k-1} \sum_{\vh \in s_q X \cap \Z^{k-1}} 1 +
\text{Error}.
$$
where
\begin{equation}
  \label{eq:17}
\text{Error} = r_q^{k-1} \sum_{\substack{ d|q \\ d>1}} \sum_{\vh
\in s_q X \cap \Z^{k-1}} e_k(\vh,d)
\end{equation}
Since $s_q = 1/r_q$, the main  term equals
$$
r_q^{k-1} \sum_{\vh \in s_q X \cap \Z^{k-1}} 1 = r_q^{k-1} \left(
\vol(s_qX) + O(s_q^{k-2})  \right) = \vol(X) + O(1/s_q).
$$
To prove the theorem we wish to show that Error$=o(1)$. To begin
with we show that the average of $e_k(\vh,d)$, over a full set of
residues modulo $d$, equals zero for $d>1$:
\begin{lem}
\label{lem:eps-average-equals-zero}
If $d>1$ then
$$
\sum_{\vh \in (\Z/d\Z)^{k-1}} e_k(\vh,d) = 0
$$
\end{lem}
\begin{proof} For any prime $p$ we have
\begin{multline*}
|\Sp|^k = \sum_{\vh \in (\Z/p\Z)^{k-1}} N_k(\vh,p) = r_p^{k-1}
|\Sp| \sum_{\vh \in (\Z/p\Z)^{k-1}} (1 + \eps_k(\vh,p)) \\
= p^{k-1} r_p^{k-1} |\Sp| + p r_p^k \sum_{\vh \in (\Z/p\Z)^{k-1}}
e_k(\vh,p)
\end{multline*}
so  that $ \sum_{\vh \in (\Z/p\Z)^{k-1}} e_k(\vh,p) = 0$. The
result follows as $e_k(\vh,d)$ is multiplicative.
\end{proof}

Throughout this section we shall take
$\tau_1=[\sqrt{2k+1/4}-\frac{1}{2}]$, $v(0)=k-2$,
$v(\tau_1)=k+\frac{1}{2}-\sqrt{2k+1/4}$, $v(\tau)=k-\tau$ for
$\tau_1+1 \leq \tau\leq k-1$ and $w(\tau)=k-9/8 + (\tau-1/2)^2/2$.

\begin{prop}
\label{prop:messy-error-bound}
Suppose that we are given
$R \in [0,1]$,  as well as $\alpha_0, 
\alpha_1, \beta_1, \alpha(\tau), \beta(\tau)>0$, for $\tau_1\leq
\tau \leq k-1$. Assume that $|\Sp|>p^{1-\alpha(\tau)}$ for
all $\tau$ and all
primes $p$ (so that $s_p\leq p^{\alpha(\tau)}$). Then 
\[
\begin{split}
\text{\rm Error} & \ll s^{\alpha_0 R-1}_q \prod_{p|q} \Big
(1+O_k(p^{1-\alpha_0}  ( A_{p,k} +   (s_p-1)/p) ) \Big )\\
& + s^{\alpha_1-\beta_1R}_q \prod_{p|q} \Big (1+ O_k(p^{\beta_1} (
A_{p,k} +  (s_p-1)/p^{1+\alpha_1} ) \Big )\\
& + \sum_{\substack{\tau=0\ \text{\rm or} \\  \tau_1\leq \tau\leq
k-1 }}
 s^{v(\tau)+\alpha(\tau) w(\tau)-(k-1)-\beta(\tau) R}_q
\prod_{p|q} \left( 1+p^{\beta(\tau)} O_k \Big( A_{p,k} + \frac{s_p
-1 }{p^{\alpha(\tau)}} \Big ) \right).
\end{split}
\]
\end{prop}

\begin{proof}
We split the divisor sum in (\ref{eq:17}) into two parts depending
on the size of the divisor $d$.

\smallskip
\noindent{\em Small $d$:}  We first consider $d \leq s_q^R$. A
point $\vh \in s_q X \cap \Z^{k-1}$ is contained in a unique
cube $C_{\vh,d} \subset \R^{k-1}$ of the form
$$
C_{\vh,d} = \{(x_1, x_2, \ldots, x_{k-1}) : dt_i \leq x_i < d(t_i+1),
t_i \in \Z, \, i =1,2, \ldots, k-1  \}
$$
We say that $\vh \in s_q X \cap \Z^{k-1}$ is a {\em $d$-interior} point
of $s_qX$ if $C_{\vh,d} \subset
s_q X$, and if $C_{\vh,d}$ intersects the boundary of $s_qX$, we say
that $h$ is a {\em $d$-boundary point} of $s_qX$.

By Lemma~\ref{lem:eps-average-equals-zero}, the sum over the
$d$-interior points is zero, and hence
\begin{equation}
  \label{eq:4}
r_q^{k-1} \sum_{\substack{ d|q \\ 1< d \leq s_q^R}} \sum_{\vh \in
s_q X \cap \Z^{k-1}} e_k(\vh,d) = r_q^{k-1} \sum_{\substack{ d|q
\\ 1< d \leq s_q^R}} \sum_{\substack{\vh \in s_q X \cap
\Z^{k-1}\\\text{$\vh$ is
      $d$-boundary point}}}
e_k(\vh,d)
\end{equation}
Now, the number of cubes $C_{\vh,d}$ intersecting the boundary of $s_qX$
is $\ll (s_q/d)^{k-2} $, and hence (\ref{eq:4}) is
$$
\ll r_q^{k-1} \sum_{\substack{ d|q \\ 1< d \leq s_q^R}}
(s_q/d)^{k-2} \sum_{\vh \in (\Z/d\Z)^{k-1}} |e_k(\vh,d)|
$$
\begin{equation}
  \label{eq:7}
= \frac{1}{s_q} \sum_{\substack{ d|q \\ 1< d \leq s_q^R}}
\frac{1}{d^{k-2}} \sum_{\vh \in (\Z/d\Z)^{k-1}} |e_k(\vh,d)|
\end{equation}

Further,
$$
\sum_{\vh \in (\Z/d\Z)^{k-1}} |e_k(\vh,d)| = \prod_{p|d} \sum_{\vh
\in (\Z/p\Z)^{k-1}} |e_k(\vh,p)|
$$

By assumption, $|e_\ell (\vh',p)| \leq A_{p,\ell} \leq A_{p,k} $
whenever $\vh'$ has $\ell\leq k$ distinct elements mod $p$.
Therefore, by (\ref{eq:106}),
\begin{equation}
  \label{eq:19}
|e_k(\vh,p)| \leq s^{k-\ell}-1 + s^{k-\ell}_p A_{p,k},
\end{equation}
for all $\vh$ with $\ell$ distinct entries modulo $p$, and so
\[
\sum_{\vh\in (\Z/p \Z)^{k-1}} |e_k(\vh,p)|  \leq
 p^{k-1} A_{p,k} + O_k \Big (\sum^{k-1}_{\ell=1} p^{k-\ell-1} (s^\ell_p-1+
 s^\ell_p A_{p,k}) \Big ) .
\]
Now $s_p/p\leq 1/2$ for $p$ large, so this error term is $\ll_k p^{k-2}
(s_p-1+s_p A_{p,k})$, and so the equation implies that
\[
\sum_{\vh \in (\Z/d \Z)^{k-1}} |e_k(h,d)|\leq
  d^{k-2} \prod_{p|d} \Big (pA_{p,k} + O_k (s_p-1+s_p A_{p,k})
  \Big).
\]
Now, $1\leq (s_q^r/d)^{\alpha_0}$ for any $\alpha_0>0$, for all
$d\leq s_q^r$, and therefore
 (\ref{eq:7}) is, for any $\alpha_0>0$,
\begin{equation}
\label{eq:y-star-two} \leq s_q^{\alpha_0 R - 1}
\prod_{p|q}\Big(1+p^{-\alpha_0}\Big(pA_{p,k}+ O_k (s_p-1+s_p
A_{p,k})\Big) \Big) ,
\end{equation}
and we get the first term in the upper bound.

\smallskip
\noindent {\em Large $d$:}  We now consider $d> s_q^R$.  Define
$\Gamma (\vh)$ as in~\ref{sec:higher-correlations-preliminaries}.
By~(\ref{eq:19}),
\[
|e_k(\vh,d)| \leq \prod_{p|d/c} A_{p,k}
\prod_{p^e\|\gamma}(s^e_p-1+s^e_p A_{p,k}),
\]
(note that $\# \{ h_0=0, h_1, \ldots, h_{k-1} \mod p\}=k - e$ if
$p|c$ but $=k$ if $p|(d/c)$), and hence 
\begin{multline*}
\sum_{\vh \in s_qX \cap \Z^{k-1}} |e_k(\vh,d)| 
\\ \leq \sum_{c|d}
(\prod_{p|d/c} A_{p,k}) \sum_{\substack{\gamma:\\
c|\gamma|c^{k-1}}} \prod_{p^e\|\gamma}(s^e_p -1+s^e_p A_{p,k})
\cdot
\sum_{\substack{\vh \in s_q X\cap \Z^{k-1} \\
\gamma(\vh)=\gamma}} 1 .
\end{multline*}
Now
$\displaystyle \sum_{\substack{\vh \in s_qX\cap \Z^{k-1}\\
    \gamma(\vh)=\gamma}} 1 \leq M_\gamma (H)
$ as defined earlier, where $H=O(s_q)$.  Using
Corollary~\ref{cor:m-gamma-bound} we bound this in various ranges:
For $\gamma \leq H$ we obtain
\begin{equation}
\label{eq:small-gamma} \ll_k H^{k-1} \sum_{c|d} (\prod_{p|d/c}
A_{p,k}) \sum_{\substack{\gamma \leq H\\ c|\gamma|c^{k-1}}}
\frac{1}{\gamma} \prod_{p^e\|\gamma} S(k,k-e) (s^e_p -1 + s^e_p
A_{p,k}) .
\end{equation}
Now, for any
$\alpha_1 >0$,
the last sum here is
\begin{multline*}
\leq \sum_{\substack{\gamma \geq 1\\ c|\gamma|c^{k-1}}} \Big
(\frac{H}{\gamma}\Big )^{\alpha_1} \frac{1}{\gamma}
\prod_{p^e\|\gamma} \left( S(k,k-e)(s^e_p -1 + s^e_p A_{p,k})
\right)
\\
= H^{\alpha_1} \prod_{p|c} \left(  \sum^{k-1}_{e=1} S(k,k-e)
\frac{s^e_p - 1+s^e_p  A_{p,k}}{p^{e(1+\alpha_1)}} \right)
\end{multline*}
and substituting this above gives that (\ref{eq:small-gamma}) is
\begin{equation}
  \label{eq:ast-3}
\ll_k H^{k-1+\alpha_1} \prod_{p|d} \Big ( A_{p,k} + O_k\Big(
\frac{s_p -1 + s_p A_{p,k}}{p^{1+\alpha_1}} \Big ) \Big )
\end{equation}
The other ranges for $\gamma$ take the form $\gamma \leq
H^{w(\tau)}$ (and $\gamma>H^{w(\tau')}$) giving a bound  $M_\gamma 
(H) \ll_k H^{v(\tau)} \prod_{p^e\|\gamma} S(k,k-e)$, and the
analogous argument then gives that the sums are, for any
$\alpha(\tau)>0$,
\begin{equation}
\label{eq:ast-4} \ll_k H^{v(\tau)+\alpha(\tau) w(\tau)}
\prod_{p|d} \Big( A_{p,k} +  O_k\Big( \frac{s_p -1 + s_p
A_{p,k}}{p^{\alpha(\tau)}} \Big )  \Big )
\end{equation}
where $\tau=0, \tau_1$ or $\tau_1+1 \leq \tau \leq k-1$.
We need to bound $r^{k-1}_q \sum_{\substack{d|q\\ d>s_q^R}} \rho(d)$
with $\rho(d)$ as in
(\ref{eq:ast-3}) or (\ref{eq:ast-4}).
Clearly this is
$$
\leq
r^{k-1}_q
\sum_{\substack{d|q\\ d\geq 1}}
\rho (d) (d/s^R_q)^{\beta}
$$ for
any $\beta >0$, and
recalling that  $H=O(s_q)$,  we obtain the bounds
\begin{equation}
\label{eq:ast-5} \ll_k s_q^{\alpha_1-\beta_1 R} \prod_{p|q} \Big
(1+p^{\beta_1} \Big ( A_{p,k} + O_k\Big( \frac{s_p -1 + s_p
A_{p,k}}{p^{1+\alpha_1}} \Big ) \Big )\Big )
\end{equation}
and 
\begin{multline}
\label{eq:ast-6} \ll_k s_q^{\nu (\tau)+\alpha(\tau) w
(\tau)-(k-1)-\beta(\tau) R} 
\cdot \\ \cdot \prod_{p|q} \left( 1+p^{\beta(\tau)}
\Big( A_{p,k} +  O_k\Big( \frac{s_p -1 + s_p
A_{p,k}}{p^{\alpha(\tau)}} \Big ) \right)
\end{multline}
for any $\alpha(\tau), \beta(\tau) >0$, where $\tau$ runs through
the relevant ranges, and the result follows.
\end{proof}

Define $\lambda_k := \min_\tau (k-1-v(\tau))/w(\tau)$ so that
$\lambda_2=(\sqrt{17}-3)/2=.56155\dots,\ \lambda_3=1/3$, and
 $\lambda_k=\frac{1}{k-1}$ for all $k\geq 4$.

We will deduce the following theorem from
Proposition~\ref{prop:messy-error-bound}, which
implies Theorem~\ref{thm:main} after the discussion in section
2.

\begin{thm}
\label{thm:main-bis} Fix $\epsilon>0$ and integer $K$. Suppose that we
are given subsets $\Sp \subset \Z/p\Z$ for each prime $p$ with
$s_p \ll_K p^{\lambda_K-\epsilon}$. Moreover assume that (1) holds
for each $k\leq K$ provided that  $0, h_1, h_2, \ldots h_{k-1}$
are distinct mod $p$.  Then, for $X \subset \{ x \in \R^{k-1} :
0<x_1 < x_2 < \ldots < x_{k-1}\}$, the $k$-level correlation
function  satisfies
$$
R_k(X, \Sq) =  \vol(X) + o_{X,k}(1)
$$
as $s_q = q/|\Sq|$ tends to infinity.
\end{thm}

This follows immediately from Proposition~\ref{prop:messy-error-bound}
and the following: 
\begin{lem}
Fix $\epsilon>0$ and assume that
\[
A_{p,k} \ll_k (1-r_p) p^{-\epsilon} \mbox{ with } s_p \ll_k
p^{\lambda_k-2\epsilon} .
\]
 Then there exists $\delta=\delta_\epsilon >0$ such
that $\text{\rm Error} \ll s^{-\delta}_q$.
\end{lem}

\begin{proof} Taking $\alpha_0=1, \alpha_1 \leq R\beta_1-2\delta$ where
$0<\beta_1<\epsilon/2$, $\beta(\tau)=0$ and
$\alpha(\tau)=\lambda_k-\epsilon$ (so that $s_p\leq
p^{\alpha(\tau)-\epsilon})$ in
Proposition~\ref{prop:messy-error-bound}, we find that the 
$p$-th term in each Euler product is $\leq
1+O((1-r_p)/p^{\epsilon/2})$.
 Now if $1\leq s_p \leq 2$ then this
is $\leq 1+O((s_p-1)/p^{\epsilon/2})=s_p^{O(1/p^{\epsilon/2})}=s_p^{o(1)}$, and
  if $s_p > 2$ this is
  $1+O(1/p^{\epsilon/2})=s_p^{O(1/p^{\epsilon/2})}=s^{o(1)}_p$.
  Thus each of the Euler products is $s_q^{o(1)}$ and the result
  follows.
\end{proof}

\section{Poisson spacings for values taken by generic polynomials}
\label{sec:general-polynomials}

Let $f$ be a polynomial of degree $n$ with integer coefficients, and
assume that $f$ has $n-1$ distinct critical values, i.e., that 
$$
\{ f(\xi) : f'(\xi) = 0, \ \xi \in \overline{\Q} \}
$$
has $n-1$ elements.  Then, for all but finitely many $p$, the set 
$$
\{ f(\xi) : f'(\xi) = 0, \ \xi \in \kbar \}
$$
also has $n-1$ elements.

We will deduce Theorem~\ref{thm:gaps-of-polynomials}
from Theorem~\ref{thm:main} together with the
following result:
\begin{thm}
\label{thm:poly-correlation}
Let $f \in \fp[x]$ be a polynomial of degree $n<p$, and let 
$$
R:= \{f(\xi) : \xi \in \kbar, f'(\xi) = 0\}.
$$  
Assume that 
$|R|=n-1$. 
If $0,h_1, h_2, \ldots h_{k-1}$ are distinct modulo 
$p$, then
$$
N_k((h_1, h_2, \ldots, h_{k-1}), p)
=
r_p^{k} \cdot p + O_{k,n}(\sqrt{p}).
$$
\end{thm}
\begin{rem}
Theorem~\ref{thm:poly-correlation} is not true for all polynomials.
For example, if we take $f(x) = x^4 - 2x^2$, then the critical values
of $f$ are $0,-1$, and for certain primes $p$, $N_2(1,p) =
3/32 \cdot p + O(\sqrt{p})$, rather than the expected answer 
$(3/8)^2 \cdot p + O(\sqrt{p})$.  See
Section~\ref{sec:counter-example} 
for more details.
\end{rem}

\subsection{Proof of Theorem~\ref{thm:poly-correlation}}
\label{sec:proof-theorem}
Assume that $n$ and $k$ are given and that $p$
is a sufficiently large prime (in terms of $n$ and $k$).
We wish to count the number of $t$ for which there exists $x_0, x_1,
\ldots x_{k-1}\in \fp$ such that
$$
f(x_i) = t+h_i \text{ for $0 \leq i \leq k-1$.}
$$
In order to study this, let $X_{k,\vh}$ be the affine curve
$$
X_{k,\vh} := \{f(x_0) = t, \ f(x_1) = t+h_1, \ldots, f(x_{k-1}) =
t+h_{k-1}    \}.
$$
and let $\fp[X_{k,\vh}]$ be the coordinate ring of 
$X_{k,\vh}$.  We then have 
\begin{multline}
\label{eq:number-of-degree-one-primes-in-K}
N_k((h_1, h_2, \ldots, h_{k-1}), p)
\\
=
|\{\m \in \fp[t] : \text{$\M|\m$ for some degree one prime $\M \in
  \fp[X_{k,\vh}]$ } \}| 
\end{multline}
In order to estimate the size of this set, we will use the Chebotarev
density theorem, made effective via the Riemann hypothesis for curves,
for the Galois closure of $\fp[X_{k,\vh}]$.  Thus, define a curve
$Y_{k,\vh}$ by letting $\fp(Y_{k,\vh})$ correspond to the Galois
closure of the extension $\fp(X_{k,\vh})/\fp(t)$.  In order to study
this extension we introduce some notation:
Given $h \in \kk$,
define a polynomial  $F_h \in \kk[x,t]$ by  
$$
F_h(x,t) := f(x) - (t+h).
$$
Since the $t$-degree of $F_h$ is one, it is 
irreducible, and thus 
$$
K_h := \kk[x,t]/F_h(x,t)
$$
is a field.
Let $L_h$ be the Galois closure of $K_h$, and let 
$$
G_h :=
\gal(L_h/\kk(t)).
$$
(Note that all field extensions considered are separable since $p>n$.)   

Hilbert has shown \cite{hilbert-full-galois-group} (e.g., see Serre
\cite{serre-topis-in-galois-theory}, chapter 4.4)  that $G_h \cong S_n$ for 
all $h$.  Our first goal is to show 
that the field extensions $L_{h_0}, \ldots, L_{h_{k-1}}$ are linearly
disjoint, or equivalently, if we let 
$$
E := L_{h_0} L_{h_1} \cdots L_{h_{k-1}}
$$
be the compositum of the fields $L_{h_0}, \ldots, L_{h_{k-1}}$, that
$\gal(E/\fp(t)) \cong S_n^k$.  

We begin with the following consequence of Goursat's Lemma:
\begin{lem}
\label{lem:subgroup-of-Sn}
Given a subset $I = \{ i_1, i_2,
\ldots, i_l\}$ of $\{1, 2, \ldots, k\}$, define a projection
  $P_I : S_n^k \to S_n^l$ by  
$$
P_I( (\sigma_1, \sigma_2, \ldots, \sigma_k) ) = 
(\sigma_{i_1}, \sigma_{i_2}, \ldots, \sigma_{i_l}).
$$
Let $K$ be a subgroup of $S_n^k$, and assume that the restriction of
$P_I$ to $K$ is surjective for all $I \subsetneq 
\{1, 2, \ldots, k\}$.  If $k>2$ then either $K = S_n^k$ or 
$$
K = \{ \sigma \in S_n^k: \sgn(\sigma) = 1
  \}.
$$
If $k=2$, there is the additional possibility that 
$$
K = \{ (\sigma_1,\sigma_2) \in S_n \times S_n : \sigma_1=\sigma_2 \},
$$
and if $k=2$ and $n=4$, we also have the  possibility that 
$$
K = \{ (\sigma_1,\sigma_2) \in S_4 \times S_4 : \sigma_1 H = \sigma_2 H \}
$$
where $H = \{1, (12)(34), (13)(24), (14)(23)\}$ is the unique
nontrivial normal subgroup of $A_4$.  
In particular, we note that if $K$ contains an odd permutation,
then $K=S_n^k$.

\end{lem}
\begin{proof}
Let $P_1=P_{\{1\}}$ be the projection on the first coordinate, put
$P_2=P_{\{2,3,\ldots k\}}$, and let 
$N_i$ be the kernel of $P_i$ restricted to $K$ for $i =1,2$.  We may
then regard $N_1$ as 
a normal subgroup of $S_n^{k-1}$, and $N_2$ as a normal subgroup of $S_n$.
By Goursat's lemma (e.g. see exercise~5 of ch.~1 in \cite{lang-algebra}),
$K$ may be described as follows (were we have identified $ S_n^k$ with
$S_n^{k-1} \times S_n$):
$$
K = \{ (x,y) \in  S_n^{k-1} \times S_n  : f_1(x) = f_2(y) \}
$$
where $f_1 : S_n^{k-1} \to S_n^{k-1}/N_1$ and $f_2 : S_n \to
S_n/N_2$ are the canonical projections, and  $S_n^{k-1}/N_1$ and
$S_n/N_2$ are identified via an isomorphism. 

We first consider the case $k>2$. 
Now, if $(\sigma_1, \sigma_2,
\ldots \sigma_{k-1}) \in N_1 \lhd S_n^{k-1}$ and $\sigma_j$ is
a transposition we find that 
$N_1 $ contains the subgroup
$$
\{ (\sigma_1, \sigma_2, \ldots, \sigma_{k-1}) : \sigma_j \in A_n
\text{ and $\sigma_i = 1$ for $i \neq j$} \}.
$$
Hence, since $P_I$ is surjective for all $I \subsetneq \{1, 2,
\ldots, k\}$, we have $A_n^{k-1} \subset N_1$.  Thus $f_1$ factors
through $S_n^{k-1}/A_n^{k-1} \cong \mathbb{F}_2^{k-1}$ and hence
$S_n^{k-1}/N_1 \cong \mathbb{F}_2^{k'} $ for some $k' < k$.  But
if $\mathbb{F}_2^{k'} \cong S_n/N_2$ then 
either $N_2 = S_n$ and $k'=0$, or $N_2 = A_n$ and $k' = 1$.  In the
first case, we find that $f_1$ and $f_2$ both are constant, and thus $K =
S_n^k$.  As for the second case, we note that
$f_2(\sigma) = \sgn(\sigma)$ and that $f_1$
must be of the form
$$
f_1((\sigma_1, \sigma_2, \ldots, \sigma_{k-1})) =
\prod_{i=1}^{k-1} \sgn(\sigma_i)^{\epsilon_i}
$$
for some choice of  $\epsilon_i \in \{0,1\}$ for $1\leq i \leq k-1$
(any homomorphism $\mathbb{F}_2^{k-1} \to \mathbb{F}_2$ is of the form
$(x_1, x_2, \ldots, x_{k-1}) \to \sum_{i=1}^{k-1} \epsilon_i x_i$).
Thus, if we put $\epsilon_k = 1$, we have
$$
K = \{(\sigma_1, \sigma_2, \ldots, \sigma_k) \in S_n^k : 
\prod_{i=1}^{k} \sgn(\sigma_i)^{\epsilon_i}  = 1\}.
$$
On the other hand, since $P_I$ is surjective for all $I \subsetneq
\{1, 2, \ldots, k\}$ we must have $\epsilon_i = 1$ for $1 \leq i \leq
k$.

As for the case $k=2$, we recall that
the only nontrivial normal subgroup of $S_n$ is $A_n$,
except when $n=4$ in which case $H$ is also a normal subgroup.  Since
$N_1$ and $N_2$ are both normal in $S_n$, and $S_n/N_1 \cong
S_n/N_2$, we must have $N_1=N_2$, and the result follows.

\end{proof}

In order to show that $\gal(E/\fp(t))$ contains an element with odd
sign, we will need the following:
\begin{lem}
\label{lem:not-small-weight}
Let $H,S \subset \fp$.  
If $p > 4^{|S|+|H|}+1$
then there exists $t \in \fp$
such that the number of $h \in H$ with $t \in S-h$ is odd.
\end{lem}
\begin{proof}
Since 
$$
|\{h \in H : t \in S-h \}| =
|\{h \in \alpha H : \alpha t \in \alpha S - h \}| 
$$
for $\alpha \in \fp^\times$, we may replace $S$ and $H$ by $\alpha
S$ and $\alpha H$ where $\alpha \in \fp^\times$ is chosen freely;
similarly we may also replace $S$ and $H$ by $S+\beta$ and $H+\beta'$
for any $\beta,\beta' \in \fp$.  Now, given $\vec{v} \in
\fp^{|S|+|H|}$ we may partition $\fp^{|S|+|H|}$ into $4^{|S|+|H|}$
boxes with sides at most $p/4$.  If $4^{|S|+|H|}<p-1$, the Dirichlet
box principle gives that there exists $\alpha', \alpha''$ such that
all components of $\alpha' \vec{v}$ and $\alpha'' \vec{v}$ differ by
at most $p/4$.  Thus, with $\alpha = \alpha' - \alpha''$ we may choose
$\beta$ such that $\alpha \vec{v} + \beta(1,1,1,\ldots, 1) \equiv
(x_1,x_2,\ldots, x_{|S|+|H|}) \mod p$ where $0 \leq x_i < p/2$ for
$1\leq i \leq |S|+|H|$.  We may thus assume that integer
representatives for all elements of $S$ can be chosen in $[0,p/2)$
and, by replacing $H$ by $H+\beta'$ for an appropriate $\beta'$, we may
also assume
that integer representatives for all elements in $H$ may be chosen in
the interval $(p/2,p]$.

Thus, if we define $h(T), s(T) \in \mathbb{F}_2[T]/(T^p-1)$ by $h(T) =
\sum_{h \in H} T^{p-h}$ and $s(T) = \sum_{s \in S} T^s$ we find that the
degrees of $h(T)$ and $s(T)$ are less than  $p/2$.
Now, if the number of $h \in H$ with $t \in S-h$ is even for all $t$,
then
$$
h(T)s(T) \equiv 0 \mod T^p-1.
$$
However, this cannot happen since the degree of $h(T)s(T)$
is less than $p$.
\end{proof}

\begin{rem}
The conclusion of the Lemma does not hold for   $p=7$, $S=
\{0,1,2,4\}$ and $H = \{0,4,6\}$, so it is necessary to make some
assumption on the size of  $p$.
\end{rem}

We can now show that the Galois group is maximal:
\begin{prop}
If $p \gg_{k,|R|} 1$ and $h_0=0, h_1, h_2, \ldots h_{k-1}$ are distinct
modulo $p$, then 
$$
\gal(E/\fp(t)) \cong S_n^k.
$$
\end{prop}

\begin{proof}
Since 
$$
\gal(E\kbar /\kbar(t))  \lhd \gal(E/\fp(t)) < S_n^k
$$
it is enough to show that $\gal(E\kbar /\kbar(t)) = S_n^k$, i.e., we
may assume that the field of constants is algebraically closed.  We
also note that this implies that the constant field of $E$ is $\fp$,
i.e., 
\begin{equation}
  \label{eq:constant-field-is-fp}
E \cap \kbar = \fp.
\end{equation}

We may regard $\gal(E\kbar /\kbar(t))$ as a subgroup of $S_n^{k-1}
\times S_n$.  By induction we may assume that the assumptions in
Lemma~\ref{lem:subgroup-of-Sn} are satisfied.  Hence $\gal(E\kbar
/\kbar(t))$ is either isomorphic to $S_n^k$, or to $\{\sigma \in S_n^k
: \sgn(\sigma) = 1\}$.  To show that the second case
cannot occur it is enough to prove that the Galois group contains
an element with odd sign.

We will now show that there exists a prime ideal $\m \subset \fp[t]$ such
that the number of $h_i$ for which $\m$ ramifies in $K_{h_i}$ is {\em
  odd}.
We begin by noting that ramification of the ideal $(t-\alpha)$ in
$K_{h_j}$ is equivalent to $\alpha + h_j \in R$.
Choose an arbitrary $r_0 \in R$.  We can then find $z \in \fp$ such
that $\m = (t-(r_0+z))$ ramifies in $K_{h_j}$ for an {\em odd} number of $j$
(for $0 \leq j \leq k-1$) in the following way:
With 
$$
R' :=  R  \cap (r_0 + \fp)
$$
we find that $ (t-(r_0+z))$ ramifies in $K_{h_j}$ if and only if
$r_0+z+h_j \in R'$.  Putting $R'' = R'-r_0$, we see that the number of
$j$ for which $r_0+z+h_j \in R'$ equals the number of $j$ for which
$z+h_j \in R''$, which in turn equals the number of $j$ such that $z
\in R''-h_j$.  By Lemma~\ref{lem:not-small-weight}, applied with
$S=R''$ and $H = \{0,h_1, \ldots, h_{k-1}\}$, it is possible to
choose $z$ so that this happens for an odd number of $j$.

If $\M$ is a prime in $E$ lying above $\m$, then the
decomposition group $\gal(E\kbar
/\kbar(t))_\M \cong \gal(E_\M/\kbar(t)_\m)$.  
After a linear change of variables we may assume the
following: $\m =(t)$,  the roots of 
$F_{h_i}(x_i,t)$
are distinct modulo $(t)$ for those $h_i$ for which 
$\m$ does not ramify in $K_{k_i}$, and for those $h_i$ for which 
$\m$ does ramify in $K_{k_i}$, we have 
$$
F_{h_i}(x_i,t) = f(x_i) - h_i - t =  x_i^2 g_i(x_i) - t
$$
where the roots of $g_i$ are distinct modulo $(t)$ and $g_i(0) \neq 0$.
Using Hensel's Lemma it readily follows that 
$E_\M = \kbar((\sqrt{t}))$, i.e., a totally ramified quadratic
extension of $\kbar(t)$.  Thus $\gal(E_\M /\kbar(t)_\m)$ is group of
order two, and is 
generated by an element $\sigma$ that maps $\sqrt{t}$ to $-\sqrt{t}$.
Now, for all $h_i$, $\sigma$ acts trivially on the unramified roots of
$F_{h_i}(x_i,t)$, and by transposing pairs of roots that are congruent
modulo $(t)$.  
Thus, when regarded as an element of $S_n^k$,  $\sigma$ is a product
of an odd number of transposition, and hence $\gal(E /\kbar(t))$ must
equal $S_n^k$. 
\end{proof}

Since $E \cap \kbar = \fp$, we note that
$$
|\{\m \in \fp[t] : \text{$\M|\m$ for some degree one prime $\M \in
  \fp[X_{k,\vh}]$ } \}|
$$
equals (taking into account  $O_{k,n}(1)$ ramified primes)
\begin{multline*}
|\{\m \in \fp[t] : \text{$\deg(\m)=1$, $\M| \m \in \fp[Y_{k,\vh}]$ and 
$\frob(\M|\m) \in \fix_{k,\vh} $}  \}|
\\+ O_{k,n}(1)
\end{multline*}
where $\fix_{k,\vh} \subset \gal(E/\fp(t))$
is the conjugacy class
\begin{multline*}
\fix_{k,\vh}  := \{ \sigma \in \gal(E/\fp(t)) \text{ such that} \\
\text{$\sigma$ fixes at least one root of $F_{h_i}$ for $i=0,
  1,\ldots, k-1$}
\}
\end{multline*}
Thus (recall Eq.~\ref{eq:number-of-degree-one-primes-in-K})
\begin{multline}
N_k((h_1, h_2, \ldots, h_{k-1}), p) 
 \\ =
|\{\m \in \fp[t] : \text{$\deg(\m)=1$, $\M| \m \in \fp[Y_{k,\vh}]$} 
\text{ and 
$\frob(\M|\m) \in \fix_{k,\vh} $}  \}|  
\\
+O_{k,n}(1)
\end{multline}
The Chebotarev density theorem
(see \cite{fried-jarden-field-arithmetic}, Proposition~5.16) gives 
$$
N_k((h_1, h_2, \ldots, h_{k-1}), p)
=
\frac{|\fix_{k,\vh}|}{|\gal(E/\fp(t))|} \cdot p + O_{k,n}(\sqrt{p}).
$$
We conclude by determining
$\frac{|\fix_{k,\vh}|}{|\gal(E/\fp(t))|}$:
\begin{lem}
If $\gal(E/\fp(t)) \cong S_n^k$
then $$
\frac{|\fix_{k,\vh}|}{|\gal(E/\fp(t))|} = 
r_p^{k} + O_{n,k}(p^{-1/2}).$$  
\end{lem}

\begin{proof}
Since $\gal(E/\fp(t)) \cong S_n^k$ we have $|\gal(E/\fp(t))| =
|S_n|^k$ and $\fix_{k,\vh}$, regarded as a subgroup of $S_n^k$, 
equals
$$
\{ (\sigma_1, \sigma_2, \ldots, \sigma_k) \in S_n^k : \sigma_i
\text{ has at least one fixed point for $1 \leq i \leq k$} \}.
$$
Thus 
$$
|\fix_{k,\vh}| =
|\{ \sigma \in S_n : \sigma
\text{ has at least one fixed point} \}|^k
$$
and hence
$$
\frac{|\fix_{k,\vh}|}{|\gal(E/\fp(t))|}
=
\left( \frac{
|\{ \sigma \in S_n : \sigma
\text{ has at least one fixed point} \}|
}{|S_n|} \right)^k
$$
Finally, again by the Riemann hypothesis for curves, we note  that
$$
r_p = |\Omega_p|/p 
$$
$$
= 
\frac{|\{t \in \kk \text{ for which there exits $x \in \kk$ such that
$f(x)=t $}\}|}{p}
$$
$$
=
\frac{|\{ \sigma \in S_n : \sigma \text{ has at least one fixed point}
  \}|}
{|S_n|}
+  O_{n,k}(p^{-1/2}).
$$
and thus 
$$
\frac{|\fix_{k,\vh}|}{|\gal(E/\fp(t))|}
 = r_p^{k} + O_{n,k}(p^{-1/2}).
$$

\end{proof}

\subsection{Theorem~\ref{thm:poly-correlation} does not hold for all
  polynomials }
\label{sec:counter-example}

We return to the example $f(x) = x^4 - 2x^2$.  The critical values of
$f$ are $0,-1$, and for $p$ large, the Galois group of the
polynomial $f(x)-t$ over $\kbar(t)$ is isomorphic to the dihedral
group $D_4$.  In fact, regarded as a subgroup of $S_4$, it is
generated by the elements $(12)(34)$ and $(23)$, corresponding to the
ramification at $t=-1$ respectively $t=0$.
However, the Galois group $H$ of the compositum of
the extensions generated by $f(x)-t$ and $f(y) - (t+1)$ is not
isomorphic to $D_4 \times D_4$; as a subgroup of $S_4 \times S_4$ it
is generated by the elements 
$(12)(34)$, $(23)(56)(78)$ and $(67)$.  This group has order $32$, and
$\fix_{2,1}$, i.e., 
the elements of $H$ that fixes at least one root of $f(x)-t$,  and at
least on root of $f(y) - (t+1)$, consists of $(), (58), (67)$. 
Thus, for primes $p$ for which the Galois group of the 
polynomials $f(x)-t$ and $f(y)-(t+1)$ over $\fp(t)$ equals the
geometric Galois group\footnote{More precisely, all sufficiently large
  primes that split completely in a certain finite extension of $\Q$,
  namely the field of constants of the
Galois extension generated by adjoining the roots of $f(x)-t$ and
$f(y)-(t+1)$ to $\Q(t)$.}, the 
following happens: 
The elements of $D_4$ that fixes at least one root of
$f(x)-t$ are $1, (14), (23)$, hence $r_p = 3/8 + O(p^{-1/2})$.  We
would thus expect that  
$$
N_2(1,p) =
r_p^2 \cdot p + O(\sqrt{p}) = 9/64 \cdot p + O(\sqrt{p}).
$$
However, since $|G'|=32$ and $|\fix_{2,1}|=3$ we have 
$$
N_2(1,p) =
3/32 \cdot p + O(\sqrt{p}).
$$

To determine for which primes $p$ splits in the field of constants (in
$\overline{\Q}$), and to determine what happens when $p$ does not
split, we ``lift'' the setup to $\Q$: Let $L_0'$ respectively $L_1'$ be the
splitting fields, over $\Q(t)$, of the polynomials $f(x) -t$
respectively $f(y)-(t+1)$.  Let $E'$ be the compositum of $L_0'$ and
$L_1'$, and let $l' = E \cap \overline{\Q}$.  Then $\gal(E'/l'(t))
\cong H$.

As before, $\gal(L_0'/(L_0' \cap \overline{\Q}) (t)) \cong D_4$ and
since it must be a normal subgroup of $S_4$, we find that $L_0' \cap
\overline{\Q} = \Q$ and that $\gal(L_0'/\Q(t)) \cong D_4$.  Similarly
$\gal(L_1'/\Q(t)) \cong D_4$, and thus $\gal(E'/\Q(t))$ embeds into
$D_4 \times D_4$, contains $H$ as a normal subgroup, hence 
$\gal(E'/\Q(t))$ is either isomorphic to $ D_4 \times D_4$ or
$H$.  We 
note that the first case is equivalent to $l'$ being a quadratic
extension of $\Q$, whereas the second is equivalent to $l' = \Q$.  On
the other hand,  $y_1 = \sqrt{1+\sqrt{t+2}}$ and $y_2 =
\sqrt{1-\sqrt{t+2}}$ are roots of $f(y) -(t+1)$, and since $\sqrt{1+t}
\in L_0' $
we find that $i \in L_0' L_1'$ since $(y_1y_2 / \sqrt{1+t})^2 =
(1-(t+2))/(1+t) = -1$.  Thus $l' = \Q(i)$ and $\gal(E'/\Q(t)) \cong D_4
\times D_4$.

Let $E$ be the splitting field of the polynomials $f(x) -t$ and 
$f(y)-(t+1)$ over $\fp$.  
Since the geometric Galois group over $\Q$ is the same as the
geometric Galois group over $\fp$ (for large $p$), reduction modulo
$p$ gives that $\gal(E/\fp(t)) \cong D_4 \times D_4$ if $p \equiv 3
\mod 4$, and $\gal(E/\fp(t)) \cong H$ if $p \equiv 1 \mod 4$ (and $p$
is sufficiently large).  Thus, as we already have seen, $N_2(1,p) =
3/32 \cdot p + O(\sqrt{p})$ if $p \equiv 1 \mod 4$.

If $p \equiv 3 \mod 4$, we have $l = E \cap \kbar = \fp(i) =
\mathbb{F}_{p^2}$, and hence the Frobenius 
automorphism must act nontrivially on $l$, i.e.,  Frobenius takes
values in
$$
\gal(E/\fp(t))^* = \{ \sigma \in \gal(E/\fp(t)) : \sigma|_l \neq 1 \}.
$$
%
Given a subset $X$ of $\gal(E/\fp(t))$, let 
\begin{multline*}
\fix(X) = \{ \sigma \in X : \text{ $\sigma$ fixes at least one root of
$f(x) = t$,} \\ \text{and at least one root of $f(y) = t+1$.} \}
\end{multline*}
The Riemann hypothesis for curves then gives that 
$$
N_2(1,p) =
\frac{|\fix(\gal(E/\fp(t))^*)|}{|\gal(E/\fp(t))^*|} \cdot p + O(\sqrt{p})
$$
Noting that $\gal(E/\mathbb{F}_{p^2}(t)) \cong H$, we conclude that
$$
|\fix(\gal(E/\fp(t))^*)| = |\fix(\gal(E/\fp(t)))| - |\fix(H)|
$$
and since $\gal(E/\fp(t) \cong D_4 \times D_4$, we find that 
$|\fix(\gal(E/\fp(t)))| = 9$.  We already know that
$|\fix(H)| = 3$, hence $|\fix(\gal(E/\fp(t))^*)| = 6$.  Moreover, since
$\gal(E/\fp(t))^* = \gal(E/\fp(t)) \setminus H $, we have
$$|\gal(E/\fp(t))^*| = |D_4 \times D_4| - |H| = 64-32 = 32,$$ and 
thus
$$
N_2(1,p) =
3/16
\cdot p + O(\sqrt{p}).
$$

In fact, this can be seen without Galois theory as follows:
Let $S_p$ be the numbers of the form $(x^2-1)^2 \mod p$.
The squares  modulo $ p$ are $b^2, 0 \leq b<p/2$, and $b^2$ is in $S_p$
iff either $(1+b)$ or $(1-b)$ is a square modulo $ p$.
Thus the number of elements of $S_p$ is (where $\jacobi{a}{p}$ is the
Legendre symbol) 
$$
\frac{1}{2} 
\sum_{b \mod p} \left( 1 - \frac{1}{4} \left( 1+\jacobi{1+b}{p} \right) 
\left( 1+\jacobi{1-b}{p} \right) \right) +O(1) 
 = \frac{3p}{8} + O(1)
$$

Now, if $a$ and $a+1$ are in $S_p$,  let $b^2=a, c^2=a+1$ so that 
$(c-b)(c+b)=1$. With $c+b=r$  we have $c=(1/2)(r+1/r)$ and $
b=(1/2)(r-1/r)$ for some value of $r \mod p$. Now 
$b^2 \in S_p$ iff either $(1/2)(2+r-1/r)$ or $(1/2)(2-r+1/r)$  is a
square modulo $ p$, and 
$c^2 \in S_p$ iff either $(1/2)(2+r+1/r)=(1/2r)(r+1)^2$ or 
$(1/2)(2-r-1/r) =(-1/2r)(r-1)^2$ is a square modulo $ p$.

On the other hand, given $r$ such that $(1/2)(2+r-1/r)$ or
$(1/2)(2-r+1/r)$  is a square modulo $p$, and $2r$ or $-2r$ is a square
modulo $ p$ then we can construct 
$a$.  (Note that $r, -r,1/r,$ and $-1/r$ lead to the same value of
$a$.)  Therefore, the  
number of $a$ such that $a$ and $a+1$ are in $S_p$ is
\begin{multline}
\label{eq:elementary}
\frac{1}{4}
\sum_{r \mod p} 
\left( 1- \frac{1}{4} \left( 1+\jacobi{2r}{p} \right) 
  \left( 1+\jacobi{-2r}{p} \right)\right)  
\cdot \\ \cdot 
\left( 1- \frac{1}{4} \left( 1+\jacobi{2r(r^2+2r-1)}{p} \right) 
  \left( 1+\jacobi{-2r(r^2-2r-1)}{p} \right) \right)
\\
= \frac{1}{64}
\sum_{r \mod p} 
\left( 9-3 \jacobi{-1}{p}+ \sum_i c_i\jacobi{f_i(r)}{p} \right)
\end{multline}
where the $f_i(r)$ are all non-constant polynomials without repeated roots
of degree $\leq 5$, and the $c_i$ are constants. By the Riemann
hypothesis for curves,  we get that $(\ref{eq:elementary})$ equals
$$
\frac{1}{64} 
  \left( 9-3 \jacobi{-1}{p} \right) p + O( p^{1/2}).
$$
Thus, if $p \equiv 1 \mod 4$ we get $N_2(1,p) = 3/32 \cdot p + O( p^{1/2})$
and if $p \equiv 3 \mod 4$ we get  $N_2(1,p) = 3/16 \cdot p + O( p^{1/2})$.

\section{Chinese Remainder Theorem for $q_1$ and $q_2$}
\label{crt-poisson-general}

By (\ref{eq:6}) we know that the spacings of elements in $\Sq$
become Poisson with parameter $\theta_q$ (as $s_q\to \infty$) if,
for any $k\geq 2$ and $X\in \mathbb B_k$, we have
$$
\sum_{\vh \in   H \cap \Z^{k-1} } \eps_k(\vh, \Sq) = o \left(
\sum_{\vh \in   H \cap \Z^{k-1} } 1\right) ,
$$
where $H=\theta_q s_q X$.
 We shall say that the spacings are {\sl strongly
Poisson}  with parameter $\theta_q$  if
$$
\sum_{\vh \in   H \cap \Z^{k-1} } \eps_k(\vh, \Sq)^2 =
o_{k}\left( \sum_{\vh \in   H \cap \Z^{k-1} } 1\right)
$$
for the same $H$.  Note that such spacings are Poisson with
parameter $\theta_q$ as may be seen by an immediate application of
the Cauchy-Schwarz inequality.

\begin{thm}
\label{thm:product-poisson-1}
 Suppose that we are given an infinite sequences of
sets $\Sa  \subset \Za$ and $\Sb \subset \Zb$ for $q_1=q_{1,n}$
and $q_2=q_{2,n}$ for all $n\geq 3$ where $(q_1,q_2)=1$. Let
$q=q_n=q_{1,n}q_{2,n}$. Suppose that the spacings of elements in
$\Sa$ become strongly Poisson with parameter $s_{q_2}$ (as $n\to
\infty$); and that
$$
\sum_{\vh \in  H \cap \Z^{k-1} } \eps_k(\vh, \Sb)^2 = O_k \left(
\sum_{\vh \in  H \cap \Z^{k-1} } 1\right)
$$
uniformly for $H\in s_q\mathbb B_k$. Then the spacing of elements in
$\Sq$ become Poisson as $n\to \infty$ if and only if the spacing
of elements in $\Sb$ become Poisson with parameter $s_{q_1}$ as
$n\to \infty$
\end{thm}
\begin{proof} By the Chinese Remainder Theorem,
\begin{multline*}
\eps_k(\vh,\Sq)+1
\\ =\frac{N_k(\vh,\Sa)}{q_1r_{q_1}^k}\frac{N_k(\vh,\Sb)}{q_2r_{q_2}^k}
 = 
\left(\eps_k(\vh,\Sa)+1\right)\left(\eps_k(\vh,\Sb)+1\right),
\end{multline*}
so that
$$
\eps_k(\vh,\Sq)=\eps_k(\vh,\Sa)\eps_k(\vh,\Sb)+\eps_k(\vh,\Sa)+\eps_k(\vh,\Sb).
$$
Now, by the Cauchy-Schwarz inequality,
\begin{multline*} 
\left| \sum_{\vh \in H \cap \Z^{k-1} }
\eps_k(\vh,\Sa)\eps_k(\vh,\Sb) \right|^2 
\\ \leq 
\left( \sum_{\vh \in H \cap \Z^{k-1} } \eps_k(\vh,\Sa)^2 \right) 
\left( \sum_{\vh \in H
\cap \Z^{k-1} } \eps_k(\vh,\Sb)^2 \right)
\\ = 
o_k\left( \left(\sum_{\vh \in  H \cap \Z^{k-1} } 1\right)^2 \right) ,
\end{multline*}
and so
$$
\sum_{\vh \in  H \cap \Z^{k-1} } \eps_k(\vh, \Sq) = \sum_{\vh \in
H \cap \Z^{k-1} } \eps_k(\vh, \Sb) + o\left( \sum_{\vh \in  H \cap
\Z^{k-1} } 1\right)
$$
by hypothesis, which gives our theorem.
\end{proof}

A simple calculation reveals that if $\Sq$ ranges over random
subsets of $\Zq$, where the probability measure on the subsets of
$\Zq$ is defined using independent Bernoulli random variables with
parameter $1/\sigma$ (see section~\ref{sec:corr-rand-select}),
then the set $\Sq$ is strongly Poisson with parameter $\theta_q>0$,
with probability 1, if and only if $\sigma=q^{o(1)}$; and thus we
can apply the above result. In fact in this case we can weaken the
hypothesis in the Theorem above:

\begin{thm} 
\label{thm:product-poisson-2}
Suppose that we are given an infinite sequences
of integers $q_1=q_{1,n}$ and $q_2=q_{2,n}$, and positive real
numbers $\sigma_1=\sigma_{q_{1,n}}, s_2=s_{q_{2,n}}$ which are
both $q_1^{o(1)}$; and let $q=q_n=q_{1,n}q_{2,n}$. We shall assume
that $\sigma_1\to \infty$ as $n\to \infty$, but not necessarily
$s_2$. Suppose $\Sb$ are given subsets of $\Zb$ with
$|\Sb|=q_2/s_2$. If $\Sa$ ranges over random subsets of $\Za$,
where the probability measure on the subsets of $\Za$ is defined
using independent Bernoulli random variables with parameter
$1/\sigma_1$ then, with probability 1, the spacing of elements in
$\Sq$ become Poisson as $n\to \infty$ if and only if the spacing
of elements in $\Sb$ become Poisson with parameter $\sigma_1$  as
$n\to \infty$.
\end{thm}

\begin{proof} The only difference from the proof above is in the
bounds we find for
$$
\left( \sum_{\vh \in  H \cap \Z^{k-1} } \eps_k(\vh,\Sa)^2 \right)
\left( \sum_{\vh \in  H \cap \Z^{k-1} } \eps_k(\vh,\Sb)^2 \right)
.
$$
Now, trivially, $N_k(\vh,\Sb) \leq N_1(0,\Sb) = |\Sb|=q_2/s_2$, and
therefore 
$|\eps_k(\vh,\Sb)| \leq s_2^{k-1}$.

If $\{ z_t: 1\leq t\leq q_1\}$ are each independent Bernoulli
random variables with parameter $1/\sigma_1$ then 
\begin{multline*}
\E( (N_k(\vh,\Sa)-q_1/\sigma_1^k)^2 ) = 
\E \left( \sum_{t \mod q_1}
\left( \prod_{i=0}^{k-1} z_{t+h_i} - \sigma_1^{-k} \right)
\right)^2 
\\ =   
\E \left( 
\sum_{t, u \mod q_1} 
\prod_{i=0}^{k-1} z_{t+h_i}
z_{u+h_i} \right) - q_1^2\sigma_1^{-2k}
\end{multline*}
Let $\eta(a)$ be the number of pairs $0\leq i,j<k$ for which
$h_j-h_i \equiv a \mod q_1$. Then
$\E \left( \sum_{t \mod q_1}
\prod_{i=0}^{k-1} z_{t+h_i} z_{t+a+h_i}\right) = q_1\sigma_1^{\eta(a)-2k}$,
so that the above equals
$$
q_1 \sigma_1^{-2k} \left( \ \sum_{a \mod q_1}
(\sigma_1^{\eta(a)}-1) \right) .
$$
Evidently $\eta(a)\leq k$ for all $a$, and there are no more than
$k^2$ values of $a$ for which $\eta(k)>0$. Thus the above is $
\ll_k  q_1 \sigma_1^{-2k}   (\sigma_1^{k}-1)$; and thus for any
$\vh \in  H$ we have $\E(\eps_k(\vh,\Sa)^2) \ll_k
\sigma_1^{k+1}/q_1$ with probability 1. The result therefore
follows since $s_2^{k-1}\sigma_1^{k+1}/q_1 = o(1)$ by hypothesis.
\end{proof}

\section{Counterexamples}
\label{sec:counterexamples}

Despite the negative aspects of Theorem~\ref{thm:product-poisson-1},
one might still hope 
that one can often take the Chinese Remainder theorem of two fairly
arbitrary sets
and obtain something that has Poisson spacings. Here we give
several examples to indicate when we cannot expect some kind of
``Central limit theorem'' for the Chinese remainder theorem!

\subsection{Counterexample 1}
\label{sec:counterexample-1}

In this case we select a vanishing proportion of the residues mod
$q_1$ randomly, together with half the residues mod $q_2$ picked
with care. Thus, in Theorem~\ref{thm:product-poisson-2} we fix $s_2=2$
and take $q_2 = 
2\SIGMA_1$ with $\Sb = \{1, 2, \ldots, \SIGMA_1\}$. Evidently
$\Sb$ is not Poisson with parameter $\sigma_1$, so $\Sq$ is not
Poisson.

\subsection{Counterexample 2}

In this case we select a vanishing proportion of the residues mod
$q_1$ and mod $q_2$ randomly, but strongly correlated. In fact, let
$u_1, u_2, \ldots, u_{q_1}$ are independent Bernoulli random
variables with probability $1/\sigma_1= q_1^{-1/2}$. Let $S=\{ i:\
u_i=1\}$, and then take $q_2 = q_1+1$ with $\Sa=\Sb=S$.

It will be convenient to let $y_i=z_i=u_i$ for $1\leq i\leq q_1$,
with $z_0=0$, and then have $y_{j+q_1}=y_j$ and $z_{j+q_2}=z_j$
for all $j$. Note that $N_2(h,\Sa) = \sum_{j=1}^{q_1} y_{j}
y_{j+h}$
 and
$N_2(h,\Sb) = \sum_{j=1}^{q_2} z_{j } z_{j+h}$ only differ by
$O(h)$ terms. (Note that  $s_2=s_1+o(1)=\sigma_1+o(1)$.)

Let $q=q_1q_2$ and define $\Sq \subset \Zq$ from $\Sa$ and
$\Sb$ using the Chinese remainder theorem, so that $j \in \Sq$ if
and only if $x_j=1$ where $x_j = y_{j} z_{j}$.

\begin{lem}
Let $I = (0,t) \subset (0,1/3)$ be an interval, and let $\Sa, \Sb$
be as above.  Then $\E( R_2(I,q)) = 2t-t^2/2 + o(1)$.
\end{lem}
\begin{proof}

Recall that
\begin{multline*}
\E( R_2(I,q)) = \sum_{h \in s_qI} \sum_{r \geq 2}^q \frac{1}{r} \E
\left( N_2(h,q) : |\Sq|=r \right) \cdot \pr(|\Sq|=r)
\end{multline*}
Since $|\Sb|= |\Sa|$ we have $|\Sq|  = |\Sa|^2$ and thus
\begin{multline*}
\E( R_2(I,q)) \\ = \sum_{h \in s_qI} \sum_{r_1=1}^{q_1}
\frac{1}{r_1^2} \E \left( \sum_{i=1}^q x_ix_{i+h} : |\Sa|=r_1
\right) \cdot \pr(|\Sa|=r_1)
\end{multline*}

Now, $\pr(|\Sa|=r_1) = (1/\sigma_1)^{r_1}(1-1/\sigma_1)^{q_1-r_1}
\binom{q_1}{r_1}$. Using the Chinese Remainder theorem and the
linearity of expectations we obtain
\begin{multline*}
\E \left( \sum_{i=1}^q x_ix_{i+h} : |\Sa|=r_1 \right) =
\sum_{i_1=1}^{q_1} \sum_{i_2=1}^{q_2} \E \left( y_{i_1} y_{i_1+h}
z_{i_2} z_{i_2+h} : |\Sa|=r_1 \right) \\
= \sum_{i_1=1}^{q_1} \sum_{i_2=1}^{q_2}
\binom{q_1-L}{r_1-L}\bigg/\binom{q_1}{r_1}
\end{multline*}
where $L = L(i_1,i_2,h)$ denotes the number of distinct integers
amongst $i_1, i_2$, the least positive residue of $i_1+h$ mod
$q_1$, and the least positive residue of $i_2+h$ mod $q_2$.
Therefore
$$
\E( R_2(I,q)) = \sum_{h \in s_qI} \sum_{i_1=1}^{q_1}
\sum_{i_2=1}^{q_2} \sum_{r_1=1}^{q_1} \frac{1}{r_1^2}
\binom{q_1-L}{r_1-L} (1/\sigma_1)^{r_1}(1-1/\sigma_1)^{q_1-r_1}.
$$
Now using, as in the proof of Lemma~\ref{lem:random-correlations}, that
\begin{multline*}
\frac 1{r_1^2} = \frac 1{(r_1-L+1)(r_1-L+2)} + \\ + O_L \left( \frac
1{(r_1-L+1)(r_1-L+2)(r_1-L+3)}  \right)
\end{multline*}
we obtain
$$
\sum_{r_1=1}^{q_1} \frac{1}{r_1^2} \binom{q_1-L}{r_1-L}
(1/\sigma_1)^{r_1}(1-1/\sigma_1)^{q_1-r_1} = \frac
1{q_1\sigma_1^L} \left( 1 + O\left( \frac 1{ \sigma_1 }
\right)\right).
$$
Moreover for each $h$ the number of $i_1,i_2$ with
$L(i_1,i_2,h)=4$ is $q_1^2+O(q_1)$, the number with $L=3$ is
$O(q_1)$, and the number with $L=2$ (which is when $i_2=i_1$) is
$q_1-h+O(1)$. Thus
\begin{multline*}
\E( R_2(I,q)) = \sum_{h \in s_qI}\left\{
\frac{q_1^2}{q_1\sigma_1^4} + \frac{O(q_1)}{q_1\sigma_1^3} +
\frac{q_1-h}{q_1\sigma_1^2} \right\} \left( 1 + O\left( \frac 1{
\sigma_1 } \right)\right)\\ = 2t-t^2/2 +O\left( \frac 1{ \sigma_1
} \right). \end{multline*}
\end{proof}

\subsection{Counterexample 3.}

In this example the sets are independently random but nonetheless,
highly correlated. We assume $m$ divides every element of $\om_1$,
a set of residues modulo $ q_1$, and every element of $\om_2$, a
set of residues modulo $ q_2$, where $m<\sigma_1, \sigma_2$ and $
\sigma_1,\sigma_2$ to be $o(\min(q_1^{1/4},q_2^{1/4}))$.

Select $x_j$'s randomly from the $q_i/m$ integers divisible by
$m$, in the range $1\le x_j\le q_i$, each selected with
probability $ m/\sigma_i$\quad ($=o(1)$, say).    Since
$N_2(h,q_i)=O(h/m) \iif m\nmid h$, and
$N_2(h,q_i)\sim|\om_i|m/\sigma_i + O(h/m)
\iif m\mid h$, we have $1+\eps_2(h,q_i)= o(1) \iif m\nmid h$, and
$1+\eps_2(h,q_i) \sim m  \iif m\mid h$.
Therefore $1+ \eps_2(h,q)=\prod_{i=1}^2 (1+ \eps_2(h,q_i)) = o(1)$
unless $m$ divides $h$, in which case it is $\sim m^2$. In
intervals (for $h$) of length $m$ this averages to $\sim \frac 1m
(m^2+o(m))=m+o(1)$ and so
$$
R_2(X,q)= 1/\sigma_q \sum_{h\in  \sigma_qX \cap \Z} \left(1 +
\eps_2(h,q) \right) \sim \frac m{\sigma_q} \vo(\sigma_q X)\sim
m\vo\ X,
$$
which is non-trivial for $m\ge 2$.

If $m_i$ divides the elements of $\om_i$, and with the elements
chosen as above then, by an analogous calculation to that above,
$$R_2(X,q)\sim \frac{m_1m_2}{\text{lcm}(m_1,m_2)} \vo(X) = \gcd(m_1,m_2)
\vo(X).
$$

\providecommand{\bysame}{\leavevmode\hbox to3em{\hrulefill}\thinspace}

\end{document}